\documentclass[12pt]{article}
\usepackage{amsmath,amsthm,amssymb,amsfonts}
\usepackage[numbers,sort&compress]{natbib}
\usepackage{enumerate}
\usepackage{color,colordvi}
\usepackage{soul}
\usepackage{fullpage}

\usepackage{changes}

\usepackage[letterpaper,top=2cm,bottom=2cm,left=3cm,right=3cm,marginparwidth=2cm]{geometry}

\usepackage{amsmath,amsfonts,amsthm,mathtools,bm}
\usepackage{graphicx}
\usepackage[colorlinks=true, allcolors=blue]{hyperref}
\usepackage[capitalise]{cleveref}
\usepackage{standalone}
\usepackage{algorithm} 
\usepackage{indentfirst}
\usepackage{enumerate}
\usepackage[affil-it]{authblk}
\usepackage{diagbox}
\usepackage{changes}
\usepackage{enumitem}
\usepackage{algpseudocode}
\usepackage{bbm}
\usepackage{orcidlink}
\usetikzlibrary{decorations.pathreplacing} % 用于绘制大括号

\theoremstyle{plain}
\newtheorem{lemma}{Lemma}[section]
\newtheorem{theorem}[lemma]{Theorem}

\newtheorem{corollary}[lemma]{Corollary}

\theoremstyle{definition}

\newtheorem{definition}[lemma]{Definition}

\newtheorem{example}[lemma]{Example}

\newtheorem{problem}[lemma]{Problem}

\newtheorem*{claim*}{Claim}

\numberwithin{equation}{section}

\newcommand{\textif}{\text{if }}
\newcommand{\textotherwise}{\text{otherwise}}

\newcommand{\RR}{\mathbb{R}}

\newcommand{\bigO}{\mathcal{O}}

\DeclareMathOperator{\dist}{\mathrm{dist}}

\DeclarePairedDelimiter{\card}{\lvert}{\rvert}
\DeclarePairedDelimiter{\set}{\lbrace}{\rbrace}

\DeclarePairedDelimiter{\paren}{\lparen}{\rparen}

\DeclarePairedDelimiter{\floor}{\lfloor}{\rfloor}
\DeclarePairedDelimiter{\ceil}{\lceil}{\rceil}

\title{Estimates of the first Dirichlet eigenvalue of graphs}

\author{Huiqiu Lin\footnote{email: huiqiulin@126.com}~\orcidlink{https://orcid.org/0000-0002-6072-4647}}
\author{Lianping Liu\footnote{email: y12252216@mail.ecust.edu.cn}~\orcidlink{0009-0007-9378-5028}}
\author{Zhe You\footnote{email: y30231280@mail.ecust.edu.cn}~\orcidlink{0009-0006-9769-5372}}
\author{Da Zhao\footnote{email: zhaoda@ecust.edu.cn}~\orcidlink{0000-0002-9582-0778}}
\affil{School of Mathematics, East China University of Science and Technology, 130 Meilong Road, Shanghai 200237, China.}

\date{}

\begin{document}

\maketitle

\begin{abstract} 
    Inspired by the Li--Yau eigenvalue-diameter estimates, we investigate lower bounds for the first Dirichlet eigenvalue in terms of the diameter (or inscribed radius) of a graph.
    Let $G = (V, E)$ be a graph with boundary $B$. 
    Assume that the interior $\Omega = V \setminus B$ is connected. 
    Let $r$ be the inscribed radius of $(G, B)$ and $d$ be the maximum degree of $G$. 
    We prove that $$\lambda_1(G, B) \geq \frac{d - 1}{r d^r},$$ which can be viewed as an analogue of the Lin–-Yau bound and the Meng–-Lin bound for normalized Dirichlet/Laplacian eigenvalues.
    We also derive the inequality $$\lambda_1(G, B) \geq \frac{1}{r |\Omega|}.$$
    In particular, for a tree $T$ with at least $3$ vertices, we show that
    $$\lambda_1(T) \geq 4 \sin^2 \frac{\pi}{4r + 6} \geq \frac{1}{(r + 1)^2}.$$
    Notably, both of the two preceding bounds are sharp up to a constant factor.
    We additionally examine upper bounds on the first Dirichlet eigenvalue under constraints on the numbers of interior and boundary vertices.
\end{abstract}



\section{Introduction}

Dirichlet eigenvalue problem constitutes a fundamental topic within spectral geometry. 
For a bounded domain $\Omega$ with smooth boundary $\partial \Omega$ in a Riemannian manifold, the Dirichlet eigenvalue problem is shown as
\begin{align}
    \Delta f + \lambda f = 0 \text{ in $\Omega$}
\end{align}
with Dirichlet boundary condition 
\begin{align}
    f|_{\partial\Omega} = 0.
\end{align}
This eigenvalue problem can also be viewed as the eigenvalue problem of the Dirichlet Laplacian operator $\Delta_\Omega$, which is defined by $\Delta_\Omega f \coloneqq (\Delta \widehat{f})|_\Omega$, where $\widehat{f}$ is an extension of $f$ to the whole manifold by assigning 0 to the region outside $\Omega$. 

Let $G = (V, E)$ be an undirected simple graph, where $V, E$ are respectively the set of vertices and edges of $G$. 
A graph with boundary is a pair $(G, B)$, where $B \subset V$ is a subset of $V$ such that $\delta\paren{V \setminus B} = B$ and $E(B, B)=\emptyset$. 
Here $\delta S \coloneqq \{x \in V \setminus S \mid x \sim y  \text{ for some } y \in S\}$ and $E (\Omega_{1}, \Omega_{2}) \coloneqq \set{\set{x, y} \in E \mid x \in \Omega_{1}, y \in \Omega_{2}}$. 
We call $B$ the boundary and $\Omega \coloneqq V \setminus B$ the interior of the graph, and assume that the subgraph on $\Omega$ is connected. 
We also write $N(x) = \set{y \in V \mid \set{x, y} \in E}$ for the set of neighbors of $x \in V$. 
The Dirichlet eigenvalue problem on graph $G$ with boundary $B$ is defined as follows.
\begin{equation*}
    \begin{cases}
        -\Delta f(x) = \lambda f(x), & x\in \Omega, \\ 
         f(x) = 0, & x\in B,
    \end{cases}
\end{equation*}
where $\Delta f(x) = \sum\limits_{y \sim x} (f(y)-f(x))$.
We call $\lambda$ the Dirichlet eigenvalue and $f$ the eigenfunction corresponding to $\lambda$.

The Dirichlet eigenvalues are ordered by $0 < \lambda_{1} < \lambda_{2} \leq \ldots \leq \lambda_{|\Omega|}$. 
For any $0 \neq f \in \RR^{\Omega}$, recall that the Rayleigh quotient is defined by
\begin{align}
    R_{(G,B)}(f) &\coloneqq \frac{\sum\limits_{(x, y) \in E(G) }(\widehat{f}(x)-\widehat{f}(y))^{2}}{\sum\limits_{x \in \Omega} f^{2}(x)},
\end{align}
where $\widehat{f}$ is obtained by extending $f$ to $\RR^V$ by $0$.

The variational characterizations of $\lambda_{k}$ are given by
\begin{align}
    \lambda_{k}(G,B) &= \min_{\substack{W \subset \mathbb{R}^{\Omega}\\ \dim W=k}} \max_{\substack{0 \neq f \in W}} R_{(G,B)}(f).
\end{align}

In the past decades, many scholars have extended classical eigenvalue problems with Dirichlet boundary condition to discrete spaces.  The Faber--Krahn inequality states that the ball has the lowest first Dirichlet eigenvalue among all bounded Euclidean domains with given volume. 
In~\cite{friedman1993some}, Friedman described the idea of a graph with
boundary, and conjectured an analogue of the Faber--Krahn inequality for regular trees, which can be regarded as discrete hyperbolic spaces. 
This conjecture was wrong and counterexamples were provided by Leydold~\cite{leydold1997faber} and Pruss~\cite{pruss1998discrete}. 
Later, Leydold \cite{leydold2002geometry} gave a complete characterization of extremal trees achieving the Faber--Krahn inequality. 
Moreover, B{\i}y{\i}ko{\u{g}}lu and Leydold~\cite{biyikouglu2007faber} considered the discrete Faber--Krahn inequality for other classes of trees with boundary. 
Recently, Bauer and Lippner~\cite{bauer2022eigenvalue} gave a discrete version of P\'olya's conjecture for the induced subgraphs of $n$-dimensional integer lattice, which is a discrete model of $\mathbb{R}^n$. 
Hua, Lin, and Su~\cite{hua2023payne} proved some analogues of Payne--P\'olya-Weinberger, Hile--Protter and Yang's inequalities for arbitrary subgraphs of $n$-dimensional integer lattice, which partially answers a question posed by Chung and Oden~\cite{chung2000weighted}. 
Furthermore, Hua and Yadin~\cite{hua2024universal} extended Hua, Lin, and Su's results to Cayley graphs of finitely generated amenable groups, as well as for the regular trees. Besides, comparisons of Dirichlet eigenvalues, Neumann eigenvalues and Laplacian eigenvalues are also important in spectral geometry~\cite{Polya1952Remarks,Payne1955Inequalities,Aviles1986Symmetry,Levine1986Levine,Friedlander1991inequalities,Ashbaugh1997Inequalities,Frank2010Inequalities,hua2024inequalities,Rohleder2025Inequalities}. 
It is also an interesting problem in the discrete settings~\cite{shi2020comparisons,Shi2025Comparison}.

For a compact Riemannian manifold $M$ of non-negative Ricci curvature with diameter $D$, Li and Yau~\cite{li1980estimates} proved that the first non-zero Laplacian eigenvalue $\mu_2$ of $M$ satisfies 
\begin{align}
    \mu_2\geq \frac{\pi^2}{2D^2},
\end{align}
and this result was further sharpened by Zhong--Yang~\cite{zhong1984estimate}. 
In the Dirichlet boundary valued problem, the diameter $D$ can be replaced by the radius of the largest geodesic ball that can be inscribed into the manifold.
The Li--Yau inequality led to many follow-up studies on the relationship between eigenvalues and the diameter. 
Let $G$ be a graph with boundary $B$ and $d$ be the (maximum) degree of $G$. 
Chung and Yau~\cite{chung2000harnack} first considered the Li--Yau type eigenvalue-diameter estimate for the first Dirichlet eigenvalue of graphs.
They proved that the first normalized Dirichlet eigenvalue $\lambda_1^\text{norm}$ of a convex subgraph $S$ in an abelian homogeneous graph $G$ satisfies 
\begin{align}
    \lambda_1^{norm} \geq \frac{1}{8 d D^2},
\end{align}
where $D$ is the diameter of $S$.
Lin--Yau~\cite{lin2010RicciCurvatureEigenvalue} proved that
\begin{align}\label{eq:lin_lower_bound}
    \mu_2^{norm} \geq \frac{1}{d D (e^{d D + 1} - 1)},
\end{align}
where $D$ is the diameter of $G$. 
Meng--Lin~\cite{meng2024LowerBoundsFirst} enhance it to 
\begin{align}\label{eq:meng_lower_bound}
    \lambda_1^{norm} \geq \frac{1}{(d + 1)^r}, \quad \text{ and } \quad \mu_2^{norm} \geq \frac{1}{(d + 1)^D},
\end{align}
where $r$ is the inscribed radius of $(G, B)$.

Motivated by their results, we focus on giving bounds of the first Dirichlet eigenvalue by inscribed radius or diameter.
The inscribed radius of $(G, B)$ is defined by $r = \sup_{v \in V} \dist(v, B)$. 
The circumscribed radius of a graph $G$ is defined by $R = \min_{v \in V} \sup_{w \in V} \dist(v,w)$. 
The diameter of a graph $G = (V, E)$ is defined by $D = \sup_{v, w \in V} \dist(v,w)$. 
Our main results are shown as follows.

\begin{theorem}\label{thm:lower_bound}
    Let $G = (V,E)$ be a graph with boundary $B$.
    Suppose $D$ is the diameter of $G$ and $r$ is the inscribed radius. 
    Let $\Omega = V \setminus B$ be the interior and assume it is connected.
    Then
    \begin{align}
        \lambda_1(G,B) \geq \frac{1}{r |\Omega|} \geq \frac{1}{D |\Omega|}. \label{Omega D-lower_bound}
    \end{align}
This bound is sharp up to a constant factor.    
\end{theorem}

We later got aware that Lenz--Stollmann~\cite{lenz2020UniversalLowerBounds} proved~\cref{Omega D-lower_bound}. 
Here we supplement that this bound is sharp up to a constant factor. 

\begin{theorem}\label{thm:lower_bound_degree_radius}
    Let $G = (V,E)$ be a graph with boundary $B$. 
    Let $d$ be the maximum degree in $G$ and $r$ be the inscribed radius of $(G, B)$. 
    Then
    \begin{align}\label{eq:lower_bound_degree_radius}
        \lambda_1(G, B) \geq \frac{d - 1}{r d^r}. 
    \end{align}
\end{theorem}

Note that the base in~\cref{eq:meng_lower_bound} is $d + 1$ while the base in~\cref{eq:lower_bound_degree_radius} is $d$. 
Also, the above two theorems are not comparable since $|\Omega|$ could be as large as $\bigO(|B| d^r)$ or as small as $r$. 

For trees, we choose all the leaves as boundary vertices. 
For convenience, we use $\lambda_1(T)$ to denote $\lambda_1(T,B)$.
We have the Li--Yau type eigenvalue estimates on trees as follows.

\begin{theorem}\label{thm:lower_bound_2}
    Let $T = (V, E)$ be a tree with leaves as boundary. 
    Let $r$ be the inscribed radius of $T$. 
    Suppose $|V| \geq 3$. 
    Then
    \begin{align}
        \lambda_1(T) \geq 4 \sin^2 \frac{\pi}{4r+2} \geq \frac{1}{r^2}  \label{l2-lower_bound_2}
    \end{align}
    if the diameter of $T$ is $2r$, and
    \begin{align}
        \lambda_1(T) \geq 4 \sin^2 \frac{\pi}{4r+6} \geq \frac{1}{(r + 1)^2}  \label{l2-lower_bound_2_odd}
    \end{align}
    if the diameter of $T$ is $2r + 1$.
    These bounds are sharp up to a constant factor.    
\end{theorem}

For trees, as the previous theorem show, we do not need any curvature condition to control the eigenvalues. 
Lin and Yau~\cite{lin2010RicciCurvatureEigenvalue} defined Bakry--\'Emery curvature on graphs. 
Hua and Lin~\cite{hua2019} classified unweighted graphs satisfying the curvature dimension condition $CD(0,\infty)$ (non-negative Bakry--\'Emery curvature) whose girth is at least five (the girth of trees are $\infty$). 
The classification only contains several graphs.
Lin, Lu, and Yau~\cite{lin2011ricci} defined another useful Ricci curvature on graphs, which is called Lin--Lu--Yau curvature.
Lin, Lu, and Yau~\cite{Lin2014ricci} provided a formula to calculate the Lin--Lu--Yau curvature of any edge $xy$ which is not in any $C_3$, $C_4$, or $C_5$, namely
\begin{align*}
    \kappa_{LLY}(x,y)=\frac{2}{d_x}+\frac{2}{d_y}-2,
\end{align*}
where $d_v$ denotes the degree of vertex $v$. 
This formula tells us that the tree with non-negative Lin--Lu--Yau curvature can only be a path, possibly infinite, or a star. 
Therefore, most of the trees do not have non-negative Ricci curvature in either sense of Bakry--\'Emery curvature or Lin--Lu--Yau curvature. 
We obtain Li--Yau type eigenvalue estimate for the first Dirichlet eigenvalue of trees with no prerequisite on non-negative Ricci curvature.

The Faber--Krahn inequality establishes the lower bound of the first Dirichlet eigenvalue. 
Next, we consider upper bound of the first Dirichlet eigenvalue for general graphs.

\begin{theorem}\label{thm:lambda1_graph}
    Let $G = (V,E)$ be a connected graph with boundary $B$ and $\Omega = V \setminus B$ be the interior.
    Then
    \begin{align}\label{eq:planar_lambda2}
        \lambda_1(G, B) \leq \frac{|E(\Omega, B)|}{\card{\Omega}}. 
    \end{align}
    The equality holds if and only if $|E(\Omega, B)|$ is divisible by $|\Omega|$ and each interior vertex is adjacent to exactly $|E(\Omega, B)|/|\Omega|$ boundary vertices.
\end{theorem}
 
As a direct corollary, we obtain a sharp upper bound for trees.
 
\begin{corollary}\label{tree-up-bound}
    Let $T$ be a tree with the set of leaves as boundary $B$. 
    Then
    \begin{align*}
        \lambda_{1}(T)\leq \frac{|B|}{|\Omega|}.
    \end{align*}   
    The equality holds if and only if $|B|$ is divisible by $|\Omega|$ and each interior vertex is adjacent to $|B|/|\Omega|$ leaves.
\end{corollary}

Note that $|\Omega| \geq D-1$, then by~\cref{tree-up-bound} we have $$\lambda_{1}\leq \frac{|B|}{|\Omega|} \leq \frac{n-D+1}{D-1}=\frac{n}{D-1}-1.$$
Then we have the following corollary.

\begin{corollary}
    Let $T$ be a tree with $n$ vertices and diameter $D$. 
    Then
    \begin{align*}
        \lambda_{1}(T) \leq \frac{n}{D-1} - 1,
    \end{align*}
    and the equality holds if and only if all interior vertices are on a unique diametral path and each interior vertex is adjacent to an equal number of leaves.
\end{corollary}

\Cref{thm:lambda1_graph} leads to a natural question: what is the best upper bound for general graphs when $|E(\Omega, B)|/|\Omega|$ is not an integer. 
Here, we consider some cases of trees in \cref{section::Upper bounds of the first Dirichlet eigenvalue}. 

The rest of this paper is organized as follows. In~\cref{section:Lower bounds of the first Dirichlet eigenvalue}, we first prove~\cref{thm:lower_bound} and~\cref{thm:lower_bound_degree_radius}.
Then, we give a stronger version of~\cref{thm:lower_bound_2}.
In~\cref{section::Upper bounds of the first Dirichlet eigenvalue}, we first prove~\cref{thm:lambda1_graph}, which is a sharp upper bound of the first Dirichlet eigenvalue for general graph.
Then, we refine this result for trees.
In~\cref{section:Concluding remark}, we present some problems.

\section{Lower bounds of the first Dirichlet eigenvalue}\label{section:Lower bounds of the first Dirichlet eigenvalue}

Let $G = (V, E)$ be a graph with boundary $B \subset V$. 
Suppose the interior $\Omega = V \setminus B$ is connected. 
Define a matrix $L_\Omega = (l_{xy})_{\Omega \times \Omega}$, where
\begin{align}
    l_{xy} = 
    \begin{cases}
        d_x, & \textif x = y, \\
        -1, & \textif x \sim y, \\
        0, &  \textotherwise.
    \end{cases}
\end{align}
Here, $d_x$ denotes the degree of vertex $x \in \Omega$ in the whole graph $G = (V, E)$ and $x \sim y$ denotes that $y$ is adjacent to $x$. 
Since $f(x) = 0$ for any $x \in B$, it is clear that any eigenvalue $\mu$ of $L_\Omega$ satisfies that 
\begin{align}
    \mu f(x) = (L_\Omega f)(x) = -\sum\limits_{y\sim x}(f(y)-f(x))= \lambda f(x), \ \forall x \in \Omega.
\end{align}
In other words, the Dirichlet eigenvalues of $(G, B)$ are eigenvalues of $L_{\Omega}$. 

We need the following lemma to determine the properties of the eigenfunctions for the first Dirichlet-Laplacian eigenvalue.

\begin{lemma}\label{simple eigenvalue-positive}
    Let $G = (V, E)$ be a graph with boundary $B$ and $\Omega = V \setminus B$ be the interior. 
    Suppose the interior is connected. 
    Let $\lambda_1(G,B)$ be the first Dirichlet eigenvalue. 
    Then the following statements hold.
    \begin{enumerate}
        \item[(i)] The eigenvalue $\lambda_1$ is a simple eigenvalue;
        \item[(ii)] there exists a positive eigenfunction corresponding to $\lambda_1$;
        \item[(iii)] all non-negative eigenfunctions of the Dirichlet Laplacian operator must belong to $\lambda_1$.
    \end{enumerate}
\end{lemma}

\begin{proof}
    Select a sufficiently large $c \geq \max_i (L_{\Omega})_{ii}$ such that $M = cI - L_{\Omega}$ is a non-negative matrix. 
    Since the subgraph on the interior is connected, $M$ is an irreducible matrix. 
    Since $L_{\Omega} = cI - M$, whenever $\mu$ is an eigenvalue of $M$, $\lambda = c - \mu$ must be an eigenvalue of $L_{\Omega}$. 
    Therefore, the smallest eigenvalue of $L_{\Omega}$ is $\lambda_1 = c - \rho$, where $\rho$ is the Perron eigenvalue of $M$, and its corresponding eigenvector is the same as the Perron eigenvector of $M$. 
    The claims follow by the Perron--Frobenius theorem.  
\end{proof}

\begin{lemma}\label{lem:partial_sum}
    Let $a_1, a_2, \ldots, a_\ell$ be a sequence of real numbers. 
    Let $s_0 = 0$, $s_i = s_{i-1} + a_i$ for $i = 1,2, \ldots, \ell$ be the partial sums. 
    Then $s_i^2 \leq i \sum_{j=1}^i a_j^2$ for $i = 1,2, \ldots, \ell$, and
    \begin{align}
        \sum_{i=1}^\ell s_i^2 \leq \frac{1}{4 \sin^2 \frac{\pi}{4\ell+2}} \sum_{i=1}^\ell a_i^2.
    \end{align}
\end{lemma}

\begin{proof}
    By Cauchy-Schwartz inequality, we have
    \begin{align}
        s_i^2 &= \paren*{\sum_{j=1}^{i} a_j}^2 \leq i \sum_{j=1}^i a_j^2.
    \end{align}
    Consider the vector $\bm{s} = (s_1, s_2, \ldots, s_\ell)^\top$ and $\bm{a} = (a_1, a_2, \ldots, a_\ell)^\top$. 
    Then $\bm{s} = L \bm{a}$, where
    \begin{align}
        L = 
        \begin{bmatrix} 
            1 & 0 & \cdots & 0 \\
            1 & 1 & \ddots & \vdots \\
            \vdots & \vdots & \ddots & 0 \\
            1 & 1 & \cdots & 1
        \end{bmatrix}.
    \end{align}
    Note that 
    \begin{align}
        \sum_{i=1}^\ell s_i^2 = \bm{a}^\top L^\top L \bm{a} \leq \lambda_{\max}(L^\top L) \bm{a}^\top \bm{a} = \lambda_{\max}(L^\top L) \sum_{j=1}^l a_j^2.
    \end{align}
    We only need to find the largest eigenvalue of $M \coloneqq L^\top L = (\min(i,j))_{1 \leq i,j \leq \ell}$. 
    Indeed, the eigenvalues of $M$ are $\frac{1}{4 \sin^2 \frac{k \pi}{4\ell+2}}$ with corresponding eigenvector $\bm{v}^{(k)} = (v_{1}^{(k)}, v_{2}^{(k)}, \ldots, \\v_{\ell}^{(k)})$ given by $v_{j}^{(k)} = \sin \frac{j k \pi}{2\ell+1}$ for $k = 1,2, \ldots, \ell$. 
\end{proof}

\subsection{Lower bound by volume and diameter}

\begin{figure}[htbp]
    \centering
    \begin{tikzpicture}[
    x=1.00mm, y=1.00mm,
    inner xsep=0pt, inner ysep=0pt,
    outer xsep=0pt, outer ysep=0pt,
    node circle/.style={line width=0.30mm, draw=L, fill=F, circle},
    node circle small/.style={line width=0.15mm, draw=L, fill=F, circle}
]

% 定义颜色
\definecolor{L}{rgb}{0,0,0}
\definecolor{F}{rgb}{0,0,0}

% 圆形节点（沿y=55水平中线对齐，保留原始坐标）
\node[node circle, minimum size=2mm, draw=L, fill=none] at (30, 55) {};
\node[node circle, minimum size=2mm] at (50, 55) {};  % 与第一个椭圆垂直对齐
\node[node circle, minimum size=2mm] at (70, 55) {};  % 与第二个椭圆垂直对齐
\node[node circle small, minimum size=1mm] at (85, 55) {};
\node[node circle small, minimum size=1mm] at (95, 55) {};
\node[node circle small, minimum size=1mm] at (105, 55) {};
\node[node circle, minimum size=2mm] at (120, 55) {}; % 与第三个椭圆垂直对齐
\node[node circle, minimum size=2mm] at (140, 55) {}; % 与第四个椭圆垂直对齐
\node[node circle, minimum size=2mm, draw=L, fill=none] at (160, 55) {};

% 椭圆节点（使用原始坐标，大小统一）
\path[line width=0.30mm, draw=L] (50,38) ellipse (6.8mm and 4.8mm);
\path[line width=0.30mm, draw=L] (70,38) ellipse (6.8mm and 4.8mm);
\path[line width=0.30mm, draw=L] (120,38) ellipse (6.8mm and 4.8mm);
\path[line width=0.30mm, draw=L] (140,38) ellipse (6.8mm and 4.8mm);

% 椭圆内标签（修正括号闭合问题）
\node at (50, 38) {$K_{\alpha-1}^{(1)}$};    
\node at (70, 38) {$K_{\alpha-1}^{(2)}$};    
\node at (120, 38) {$K_{\alpha-1}^{(2\ell-2)}$};
\node at (140, 38) {$K_{\alpha-1}^{(2\ell-1)}$};

% 端点标签（保留原始位置）
\node at (30, 58) {$v_0$};
\node at (50, 58) {$v_1$};
\node at (70, 58) {$v_2$};
\node at (120, 58) {$v_{2\ell-2}$};
\node at (140, 58) {$v_{2\ell-1}$};
\node at (160, 58) {$v_{2\ell}$};

% 路径连接（保留原始规整化路径）
\path[line width=0.30mm, draw=L] (31,55) -- (75,55);
\path[line width=0.30mm, draw=L] (115,55) -- (159,55);

\path[line width=0.30mm, draw=L] (50,55) -- (50-4, 38+3);  % 椭圆1左线
\path[line width=0.30mm, draw=L] (50,55) -- (50+4, 38+3);  % 椭圆1右线

\path[line width=0.30mm, draw=L] (70,55) -- (70-4, 38+3);  % 椭圆2左线
\path[line width=0.30mm, draw=L] (70,55) -- (70+4, 38+3);  % 椭圆2右线

\path[line width=0.30mm, draw=L] (120,55) -- (120-4, 38+3);% 椭圆3左线
\path[line width=0.30mm, draw=L] (120,55) -- (120+4, 38+3);% 椭圆3右线

\path[line width=0.30mm, draw=L] (140,55) -- (140-4, 38+3);% 椭圆4左线
\path[line width=0.30mm, draw=L] (140,55) -- (140+4, 38+3);% 椭圆4右线

\end{tikzpicture}
    \caption{The path-cliques graph $PC(\ell, \alpha)$}
    \label{fig:a path-cliques graph}
\end{figure}
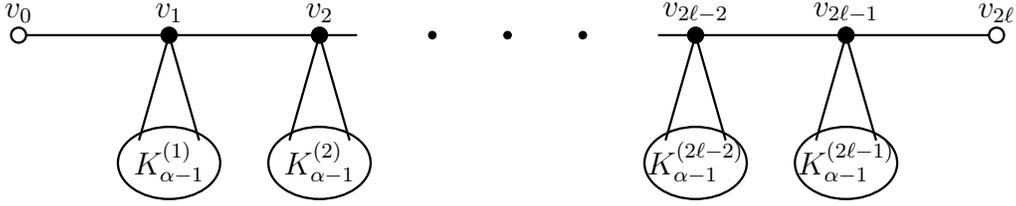

\begin{definition}\label{def:Path-Cliques Graph}
    Let $\ell \geq 2$ and $\alpha \geq 1$ be integers.
    The path-cliques graph $PC(\ell, \alpha)$ is obtained through the following procedure.
    Let $P_{2\ell}$ be the path $v_0 \sim v_1 \sim \cdots \sim v_{2\ell}$ of length $2\ell$.
    For each vertex $v_i, 1 \leq i \leq \ell$, connect $v_i$ to every vertex in a copy of $(\alpha-1)$-clique $K_{\alpha-1}^{(i)}$.
    See~\cref{fig:a path-cliques graph}.
\end{definition}

% The following example shows that \replaced{the lower bound in~\cref{Omega D-lower_bound} is sharp up to a constant factor}{, in a certain sense, Inequality \ref{Omega D-lower_bound} is as good as possible}.

\begin{lemma}\label{lemma::sharp up to constant for general graph}
    Let $G = (V, E)$ be the path-cliques graph $PC(\ell, \alpha)$. 
    Let the set of pendent vertices be the boundary $B$ of $G$ and $\Omega = V \setminus B$ be the interior. 
    Let $D$ be the diameter of $G$. 
    Then $\lambda_{1}(G, B) \leq \frac{12}{|\Omega|D}$, where $|\Omega| = (2\ell - 1)\alpha$ and $D = 2\ell$.
\end{lemma}

\begin{proof}
    Let $U_j \coloneqq \set{v_j} \cup V(K_{\alpha-1}^{(j)})$ for $1 \leq j \leq 2\ell - 1$.
    Define a test function $f$ by $f(v) = j$ for any $v \in U_j, U_{2\ell-j} \ (1 \leq j \leq \ell)$.
    Note that $f$ vanishes on the boundary. 
    Then
    \begin{align*}
        \lambda_1(G, B) &\leq R_{(G, B)}(f) \\
        &= \frac{2\ell}{\alpha\ell^2 + 2\alpha\sum_{j=1}^{\ell-1} j^2}\\
        &= \frac{6}{3\alpha\ell + \alpha(\ell - 1)(2\ell - 1)}\\
        &< \frac{6}{(\ell-1)|\Omega|} < \frac{12}{|\Omega|D}. \qedhere
    \end{align*}
\end{proof}

Now we prove our first main result. 

\begin{proof}[\bf{{Proof of~\cref{thm:lower_bound}}}]
    Let $f \in \RR^\Omega$ be a positive eigenfunction of the first Dirichlet eigenvalue by Lemma \ref{simple eigenvalue-positive} and $\widehat{f}$ is obtained by extending $f$ to $\RR^V$ by $0$.
    Suppose that $f$ achieves its maximum at a vertex $v^* \in \Omega$ and $\dist(v^*, B) = \ell$. 
    Consider a shortest path $P = (V_P, E_P)$ from $v^*$ to $v_0 \in B$ of length $\ell \leq r$, say $v^* = v_\ell \sim v_{\ell-1} \sim \cdots \sim v_0$.
    Then $a_i = \widehat{f}(v_i) - \widehat{f}(v_{i-1})$ is a real sequence of length $\ell$. 
    By~\cref{lem:partial_sum} we have $\sum_{(x, y) \in E_P} (\widehat{f}(x) - \widehat{f}(y))^2 \geq \frac{1}{\ell} f^2(v^*)$.
    Hence,
    \begin{align}
        \lambda_1 &= \frac{\sum_{(x,y) \in E} (\widehat{f}(x) - \widehat{f}(y))^2}{\sum_{x \in \Omega} f^2(x)}\\
        &\geq \frac{\sum_{(x,y) \in E_P} (\widehat{f}(x) - \widehat{f}(y))^2}{|\Omega| f^2(v^*)}\\
        &\geq \frac{\frac{1}{\ell} f^2(v^*)}{|\Omega| f^2(v^*)} \\
        &= \frac{1}{|\Omega| \ell} \geq \frac{1}{|\Omega|r} \geq \frac{1}{|\Omega| D}. 
    \end{align}
    This bound is sharp up to a constant factor due to~\cref{lemma::sharp up to constant for general graph}.
\end{proof}

\subsection{Lower bound by degree and diameter}

Let $G = (V, E)$ be a graph with boundary $B$. 
Let $\Omega = V \setminus B$ be the interior. 
An interior-boundary path in a graph $G$ with boundary is a path connecting an interior vertex and a boundary vertex. 
Let $\mathcal{P}$ be a collection of paths in $G$, where each path $P \in \mathcal{P}$ is a pair $P = (V_P, E_P)$. 
We say $\mathcal{P}$ is a $c$-covering on vertices if each interior vertex $v \in \Omega$ occurs at least $c$ times in the union of sets $V_P, P \in \mathcal{P}$ as multi-sets. 
We say $\mathcal{P}$ is a $p$-packing on edges if each edge $e \in E$ occurs at most $p$ times in the union of $E_P, P \in \mathcal{P}$ as multi-sets. 

\begin{lemma}\label{lem:graph_edge_packing}
    Let $G = (V, E)$ be a graph with boundary $B$. 
    Let $d \geq 2$ be the maximum degree in $G$ and $r$ be the inscribed radius of $(G, B)$. 
    Then $G$ has a collection $\mathcal{P}$ of interior-boundary paths such that $\mathcal{P}$ is a $1$-covering on vertices and $\frac{d^r - 1}{d - 1}$-packing on edges with each path $P \in \mathcal{P}$ of length at most $r$.
\end{lemma}

\begin{proof}
    Let $\Omega = V \setminus B$ be the interior. 
    For every $v \in \Omega$, we have $\dist(v, B) \leq r$. 
    Therefore, there exists one shortest interior-boundary path $P_v : v = v_0 \sim v_1 \sim \cdots \sim v_{\ell} \in B$ of length $\ell = \dist(v, B)$. 
    Denote by $e_v$ the arc $v_0 \rightarrow v_1$. 
    We construct a directed spanning forest $T = (V, A)$ with $A$ the set of arcs $e_v$, $v \in \Omega$. 
    For each $v \in \Omega$, let $P_v'$ be the shortest interior-boundary path in $T$ from $v$. 
    Let $\mathcal{P}$ be the collection of interior-boundary paths $P_v'$, $v \in \Omega$. 
    It is clear that $\mathcal{P}$ is a $1$-covering on vertices in $(G, B)$ and that each $P_v' \in \mathcal{P}$ is of length at most $r$. 
    We claim that it is also a $\frac{d^r - 1}{d - 1}$-packing on edges in $G$. 
    Let $c_e$ be the number of occurrences of an edge $e \in E$ in $\mathcal{P}$. 
    If an edge $e$ does not appear as an arc in $T$, it is clear that $c_e = 0 < d^r$. 
    If an edge $e = vw$ appears as an arc $v \rightarrow w$ in $T$, then $c_e$ is the number of ancestors (including itself) of $v$ in $T$. 
    Since the maximum degree in $G$ is $d$, the length of each path $P \in \mathcal{P}$ is at most $r$, we have that 
    \begin{equation*}
        c_e \leq 1 + (d - 1) + \cdots + (d - 1)^{r-1} \leq \frac{d^r - 1}{d - 1}. \qedhere
    \end{equation*}
\end{proof}

Now we prove the second main theorem.

\begin{proof}[\bf{{Proof of~\cref{thm:lower_bound_degree_radius}}}]
    By~\cref{lem:graph_edge_packing}, we have a collection $\mathcal{P}$ of interior-boundary paths such that $\mathcal{P}$ is a $1$-covering on vertices and $\frac{d^r-1}{d-1}$-packing on edges with each path $P \in \mathcal{P}$ of length at most $r$.
    Let $f$ be an eigenfunction of the first Dirichlet eigenvalue.
    Note that by~\cref{lem:partial_sum} for each path $P = (V_P, E_P) \in \mathcal{P}$, we have $f^2(v) \leq r \sum_{(x, y) \in E_P} (f(x) - f(y))^2$ for $v \in V_P$. 
    For each edge $e = (x, y) \in E$, the number of occurrences of $e$ in $\mathcal{P}$ is $c_e \leq \frac{d^r-1}{d-1}$. 
    Hence
    \begin{align*}
        \lambda_{1} &= \frac{\sum_{(x,y) \in E} (f(x) - f(y))^{2}}{\sum_{x \in \Omega} f^{2}(x)} \\
        &\geq \frac{\sum_{(x,y) \in E} (f(x) - f(y))^{2}}{\sum_{P \in \mathcal{P}} \sum_{x \in V_P} f^2(x)} \\
        &\geq \frac{\sum_{(x,y) \in E} (f(x) - f(y))^{2}}{r \sum_{e=(x, y) \in E} c_e (f(x) - f(y))^{2}} \\
        &\geq \frac{\sum_{(x,y) \in E} (f(x) - f(y))^{2}}{r \frac{d^r-1}{d-1} \sum_{(x, y) \in E} (f(x) - f(y))^{2}} \\
        &\geq \frac{d - 1}{r d^r} . \qedhere 
    \end{align*}
\end{proof}

Given a tree $T$ with inscribed radius $r$, a center of the tree $T$ is a vertex $v$ whose distance to the leaf set is exactly $r$. 
For a tree, it either has exactly one center or exactly two centers connected by an edge~\cite[Theorem 4.2]{harary2001GraphTheory}. 

\begin{lemma}\label{lem:tree_decomposition}
    Let $T = (V, E)$ be a tree of at least three vertices with leaves as boundary. 
    Let $r$ be the inscribed radius of $T$. 
    Then $T$ has a collection $\mathcal{P}$ of interior-boundary paths such that $\mathcal{P}$ is a $1$-covering on vertices and $1$-packing on edges with each path $P_i$ of length at most $r$ if $T$ has exactly one center; and $T$ has a collection $\mathcal{P}$ of interior-boundary paths such that $\mathcal{P}$ is a $1$-covering on vertices and $1$-packing on edges with each path $P_i$ of length at most $r+1$ if $T$ has exactly two centers.
\end{lemma}

\begin{proof}
    A rooted tree with boundary, denoted by $(T, v, B)$, is a tree $T = (V, E)$ together with a vertex $v \in V$ and $B$ a subset of leaves. 
    We assign the tree $T$ with a root by one of its centers $v$, namely a vertex whose distance to the leaf set is exactly the inscribed radius $a$, and boundary $B$ the set of leaves. 
    Since $T$ has at least $3$ vertices, its interior is non-empty. 
    Hence, $v$ is an interior vertex. 
    \begin{enumerate}
        \item Define the initial collection of rooted trees with boundary as $\mathcal{T} = \set{(T, v, B)}$ and the initial collection of interior-boundary paths as $\mathcal{P} = \emptyset$.
        \item If $\mathcal{T}$ is non-empty, then take one rooted tree with boundary from $\mathcal{T}$, say $(T', v', B')$ such that the root $v'$ is closest to the center $v$ among all rooted trees in $\mathcal{T}$.  
            Take one shortest path $P$ from $v'$ to one of the boundary vertex in $B'$. 
            Add $P$ to $\mathcal{P}$. 
            Remove the edges in $P$ from $T'$, and we obtain a forest $F$. 
            For each maximal tree $T''$ (with at least one edge) in the forest $F$, we assign a root by connecting vertex $v''$ on the path $P$ and boundary $B'' = B' \cap V_{T''}$. 
            Delete $(T', v', B')$ from $\mathcal{T}$ and add those new rooted trees with boundary $(T'', v'', B'')$ to $\mathcal{T}$. 
        \item Repeat the above step until $\mathcal{T}$ is empty. 
    \end{enumerate}
    For each $(T', v', B')$ appeared in the procedure, the vertex $v'$ must be an interior vertex of $T$, because no rooted tree in the forest grows on a boundary vertex. 
    The procedure will terminate because in each step, the total number of edges in $\mathcal{T}$ decreases. 
    If $T$ has exactly one center $c$, then for each path $P \in \mathcal{P}$ it may start with but can not cross the center. 
    Hence, the length of $P$ is at most $a$.
    If $T$ has exactly two centers $c_1, c_2$ with connecting edge $c_1 \sim c_2$, then for each path $P \in \mathcal{P}$, either it does not cross the edge $c_1 \sim c_2$ or does cross the edge $c_1 \sim c_2$. 
    If $P$ does not cross the edge $c_1 \sim c_2$, then the length of $P$ is at most $r$.
    If $P$ does cross the edge $c_1 \sim c_2$, then the length of $P$ can only be $r+1$ as it can only start with the edge $c_1 \sim c_2$ (or $c_2 \sim c_1$). 
\end{proof}

Then we will prove a strengthened version of~\cref{thm:lower_bound_2}. 

\begin{theorem}\label{thm:lower_bound_2_strength}
    Let $G = (V,E)$ be a graph with boundary $B$. 
    Let $\mathcal{P}$ be a collection of interior-boundary paths.
    Suppose $\mathcal{P}$ is a $c$-covering on vertices and $p$-packing on edges with each path $P \in \mathcal{P}$ of length at most $\ell$. 
    Then
    \begin{align}
        \lambda_1(G,B) \geq \frac{4c}{p} \sin^2 \frac{\pi}{4\ell+2} \geq \frac{c}{p \ell^2}. 
    \end{align}
\end{theorem}

\begin{proof}[{{Proof of~\cref{thm:lower_bound_2_strength}}}]
    Let $f$ be an eigenfunction of the first Dirichlet eigenvalue.
    Note that by~\cref{lem:partial_sum} for each path $P = (V_P, E_P) \in \mathcal{P}$, we have 
    \begin{align}
        \sum_{v \in V_P} f^2(v) \leq \frac{1}{4 \sin^2 \frac{\pi}{4|E_P|+2}}  \sum_{(x, y) \in E_P} (f(x) - f(y))^2.
    \end{align}
    Hence
    \begin{align*}
        \lambda_{1} &= \frac{\sum_{(x,y) \in E} (f(x) - f(y))^{2}}{\sum_{x \in \Omega} f^{2}(x)} \\
        &\geq \frac{\sum_{(x,y) \in E} (f(x) - f(y))^{2}}{\frac{1}{c}\sum_{P \in \mathcal{P}} \sum_{v \in V_P} f^{2}(v)}\\
        &\geq \frac{c \sum_{(x,y) \in E} (f(x) - f(y))^{2}}{\sum_{P \in \mathcal{P}} \frac{1}{4 \sin^2 \frac{\pi}{4|E_P|+2}} \sum_{(x, y) \in E_P} (f(x) - f(y))^{2}} \\
        &\geq \frac{c \sum_{(x,y) \in E} (f(x) - f(y))^{2}}{ \frac{1}{4 \sin^2 \frac{\pi}{4\ell+2}} \sum_{P \in \mathcal{P}} \sum_{(x, y) \in E_P} (f(x) - f(y))^{2}} \\
        &\geq \frac{4c \sum_{(x,y) \in E} (f(x) - f(y))^{2}}{p \sum_{(x,y) \in E} (f(x) - f(y))^{2}} \sin^2 \frac{\pi}{4\ell+2} \\
        &= \frac{4c}{p} \sin^2 \frac{\pi}{4\ell+2}. \qedhere 
    \end{align*}
\end{proof}

Now we prove the third main theorem.

\begin{proof}[\bf{{Proof of~\cref{thm:lower_bound_2}}}]
    We apply~\cref{lem:tree_decomposition,thm:lower_bound_2_strength}. 
\end{proof}

\section{Upper bounds of the first Dirichlet eigenvalue}\label{section::Upper bounds of the first Dirichlet eigenvalue}

In this section, we denote the neighbor set of any vertex $v\in V$ by $N(v)$.

\subsection{Upper bound by volume}

\begin{lemma}[{\cite[Theorem 4.2.2]{Horn_Johnson_2012}}]\label{lem:extreme_eigenfunction}
    Let $M$ be a $n \times n$ real symmetric matrix. 
    Let $\lambda_1 < \lambda_2 < \ldots < \lambda_k$ be the distinct eigenvalues of $M$. 
    If $x^\top M x = \lambda_1 x^\top x$ for $x \in \RR^n$, then $Mx = \lambda_1 x$. 
    If $x^\top M x = \lambda_k x^\top x$ for $x \in \RR^n$, then $Mx = \lambda_k x$. 
\end{lemma}

\begin{lemma}\label{lambda at most minimum degree}
    Let $G$ be a connected graph with boundary $B$. 
    Let $\Omega = V \setminus B$ be the interior. 
    Let $\delta \geq 2$ be the minimum degree among interior vertices.
    Suppose $k = |\Omega| \geq 1$.
    Then $\lambda_1(G,B) \leq \delta$.
    The equality holds if and only if $G$ is a star graph.
\end{lemma}

\begin{proof}
    Let $v_0$ be an interior vertex of degree $\delta$.
    Define a test function as follows.
    \begin{align*}
        f(v) = 
        \begin{cases}
            1, & \textif v = v_0, \\
            0, &  \textotherwise.
        \end{cases}
    \end{align*}
    Note that $f$ vanishes on the boundary. 
    Then
    \begin{align*}
        \lambda_1(G,B)\leq R_{(G,B)}(f)=\delta.
    \end{align*}
    By~\cref{simple eigenvalue-positive,lem:extreme_eigenfunction}, the equality holds if and only if $G$ is a star.
\end{proof}

Now we prove the fourth main theorem.

\begin{proof}[\bf{{Proof of Theorem \ref{thm:lambda1_graph}}}]
    Choose $1_{\Omega}$ as a test function, which denotes the indicator function that is equal to $1$ on all the interior vertices. 
    Then 
    \begin{align}
        \lambda_{1}(G,B) \leq \frac{|E(\Omega, B)|}{|\Omega|}.
    \end{align}

    If the equality holds, then by~\cref{lem:extreme_eigenfunction} we have that $1_{\Omega}$ is a $\lambda_1$-eigenfunction, where $\lambda_1 = |E(\Omega, B)|/|\Omega|$. 

    Based on the definition of Dirichlet eigenvalues, we have
    \begin{align}
        \frac{|E(\Omega, B)|}{|\Omega|} \cdot 1_{\Omega}(x) = - \Delta_{\Omega} 1_{\Omega}(x)
    \end{align}
    for any $x \in \Omega$, implying that each interior vertex is adjacent to exactly $|E(\Omega, B)|/|\Omega|$ boundary vertices.
    Therefore, the equality holds if and only if $|E(\Omega, B)|/|\Omega|$ is an integer and each interior vertex is adjacent to exactly $|E(\Omega, B)|/|\Omega|$ boundary vertices.
\end{proof}

\subsection{Upper bound for tree}

Let $T$ be a tree with $b \geq 2$ leaves and $k \geq 1$ interior vertices. 
If $b/k$ is an integer, then~\cref{tree-up-bound} characterizes the extremal trees. 
In the following we investigate the extremal trees where $b/k$ is not an integer. 
Note that by~\cref{tree-up-bound} we have $\lambda_1(T) \leq \frac{b}{k} < a + 1$, where $a \coloneqq \floor{\frac{b}{k}}$.

In the following we consider a family of trees called star-like path tree. 

\begin{figure}[htbp]
    \centering
    \begin{tikzpicture}[x=2.00mm, y=2.00mm, inner xsep=0pt, inner ysep=0pt, outer xsep=0pt, outer ysep=0pt, scale=0.8]
    % \path[line width=0mm] (8.14,37.50) rectangle +(142.41,56.89);  % 横坐标+3（原5.14→8.14）
    \definecolor{L}{rgb}{0,0,0}
    % \path[line width=0.30mm, draw=L] (11,92.39);  % 横坐标+3（原8→11）
    \definecolor{F}{rgb}{0,0,0}
    \path[line width=0.30mm, draw=L, fill=F] (38,65) circle (0.50mm) node[below=4pt] {$u_0$};  % 35→38
    \path[line width=0.30mm, draw=L, fill=F] (48,65) circle (0.50mm) node[below=4pt] {$u_1$};  % 45→48
    \path[line width=0.30mm, draw=L, fill=F] (53,65) circle (0.20mm);  % 50→53
    \path[line width=0.30mm, draw=L, fill=F] (58,65) circle (0.20mm);  % 55→58
    \path[line width=0.30mm, draw=L, fill=F] (63,65) circle (0.20mm);  % 60→63
    \path[line width=0.30mm, draw=L, fill=F] (68,65) circle (0.50mm) node[below=4pt] {$u_{c-1}$};  % 65→68
    \path[line width=0.30mm, draw=L, fill=F] (78,65) circle (0.50mm) node[below=4pt] {$u_c$};  % 75→78
    \path[line width=0.30mm, draw=L] (68,65) -- (78,65);  % 65→68，75→78
    \path[line width=0.30mm, draw=L, fill=F] (88,80) circle (0.50mm) node[left=4pt] {$u_{c,1}$};  % 85→88
    \path[line width=0.30mm, draw=L, fill=F] (88,50) circle (0.50mm) node[left=4pt] {$u_{c,e}$};  % 85→88
    
    \draw[decorate,decoration={brace,amplitude=5pt,mirror}] (89,52) -- (89,78) node[midway,xshift=10pt] {$e$};  % 86→89
     
    \path[line width=0.30mm, draw=L] (93,82) circle (0.50mm);  % 90→93
    \path[line width=0.30mm, draw=L, fill=F] (93,81) circle (0.20mm);  % 90→93
    \path[line width=0.30mm, draw=L, fill=F] (93,80) circle (0.20mm);  % 90→93
    \path[line width=0.30mm, draw=L, fill=F] (93,79) circle (0.20mm);  % 90→93
    \path[line width=0.30mm, draw=L] (93,78) circle (0.50mm);  % 90→93
    \path[line width=0.30mm, draw=L] (93,48) circle (0.50mm);  % 90→93
    
    \draw[decorate,decoration={brace,amplitude=5pt}] (93.5,82.5) -- (93.5,77.5) node[midway,xshift=10pt] {$q$};  % 90.5→93.5
    \draw[decorate,decoration={brace,amplitude=5pt}] (93.5,52.5) -- (93.5,47.5) node[midway,xshift=10pt] {$q$};  % 90.5→93.5
    
    \path[line width=0.30mm, draw=L, fill=F] (93,49) circle (0.20mm);  % 90→93
    \path[line width=0.30mm, draw=L, fill=F] (93,50) circle (0.20mm);  % 90→93
    \path[line width=0.30mm, draw=L, fill=F] (93,51) circle (0.20mm);  % 90→93
    \path[line width=0.30mm, draw=L] (93,52) circle (0.50mm);  % 90→93
    \path[line width=0.30mm, draw=L, fill=F] (88,75) circle (0.20mm);  % 85→88
    \path[line width=0.30mm, draw=L, fill=F] (88,71) circle (0.20mm);  % 85→88
    \path[line width=0.30mm, draw=L, fill=F] (88,67) circle (0.20mm);  % 85→88
    \path[line width=0.30mm, draw=L, fill=F] (88,63) circle (0.20mm);  % 85→88
    \path[line width=0.30mm, draw=L, fill=F] (88,59) circle (0.20mm);  % 85→88
    \path[line width=0.30mm, draw=L, fill=F] (88,55) circle (0.20mm);  % 85→88
    \path[line width=0.30mm, draw=L] (88,80) -- (93,82);  % 85→88，90→93
    \path[line width=0.30mm, draw=L] (88,80) -- (93,78);  % 85→88，90→93
    \path[line width=0.30mm, draw=L] (88,50) -- (93,48);  % 85→88，90→93
    \path[line width=0.30mm, draw=L] (88,50) -- (93,52);  % 85→88，90→93
    \path[line width=0.30mm, draw=L] (88,80) -- (78,65);  % 85→88，75→78
    \path[line width=0.30mm, draw=L] (88,50) -- (78,65);  % 85→88，75→78
    \path[line width=0.30mm, draw=L, fill=F] (28,80) circle (0.50mm) node[right=4pt] {$u_{0,1}$};  % 25→28
    \draw[decorate,decoration={brace,amplitude=5pt, mirror}] (22.5,82.5) -- (22.5,77.5) node[midway,xshift=-10pt] {$q$};  % 19.5→22.5
    \draw[decorate,decoration={brace,amplitude=5pt, mirror}] (22.5,52.5) -- (22.5,47.5) node[midway,xshift=-10pt] {$q$};  % 19.5→22.5
    
    \path[line width=0.30mm, draw=L, fill=F] (28,75) circle (0.20mm);  % 25→28
    \path[line width=0.30mm, draw=L, fill=F] (28,71) circle (0.20mm);  % 25→28
    \path[line width=0.30mm, draw=L, fill=F] (28,67) circle (0.20mm);  % 25→28
    \path[line width=0.30mm, draw=L, fill=F] (28,63) circle (0.20mm);  % 25→28
    \path[line width=0.30mm, draw=L, fill=F] (28,59) circle (0.20mm);  % 25→28
    \path[line width=0.30mm, draw=L, fill=F] (28,55) circle (0.20mm);  % 25→28
    \path[line width=0.30mm, draw=L] (23,82) circle (0.50mm);  % 20→23
    \path[line width=0.30mm, draw=L, fill=F] (23,81) circle (0.20mm);  % 20→23
    \path[line width=0.30mm, draw=L, fill=F] (23,80) circle (0.20mm);  % 20→23
    \path[line width=0.30mm, draw=L, fill=F] (23,79) circle (0.20mm);  % 20→23
    \path[line width=0.30mm, draw=L] (23,78) circle (0.50mm);  % 20→23
    \path[line width=0.30mm, draw=L, fill=F] (28,50) circle (0.50mm) node[right=4pt] {$u_{0,d}$};  % 25→28
    \draw[decorate,decoration={brace,amplitude=5pt}] (27,52) -- (27,78) node[midway,xshift=-10pt] {$d$};  % 24→27
    
    \path[line width=0.30mm, draw=L] (23,52) circle (0.50mm);  % 20→23
    \path[line width=0.30mm, draw=L, fill=F] (23,49) circle (0.20mm);  % 20→23
    \path[line width=0.30mm, draw=L, fill=F] (23,50) circle (0.20mm);  % 20→23
    \path[line width=0.30mm, draw=L, fill=F] (23,51) circle (0.20mm);  % 20→23
    \path[line width=0.30mm, draw=L] (23,48) circle (0.50mm);  % 20→23
    \path[line width=0.30mm, draw=L] (36,70) circle (0.50mm);  % 33→36
    \path[line width=0.30mm, draw=L] (38,65) -- (36,70);  % 35→38，33→36
    \path[line width=0.30mm, draw=L, fill=F] (37,70) circle (0.20mm);  % 34→37
    \path[line width=0.30mm, draw=L, fill=F] (38,70) circle (0.20mm);  % 35→38
    \path[line width=0.30mm, draw=L, fill=F] (39,70) circle (0.20mm);  % 36→39
    \path[line width=0.30mm, draw=L] (40,70) circle (0.50mm);  % 37→40
    \draw[decorate,decoration={brace,amplitude=5pt}] (35.5,71) -- (40.5,71) node[midway,yshift=10pt] {$p$};  % 32.5→35.5，37.5→40.5
    \draw[decorate,decoration={brace,amplitude=5pt}] (45.5,71) -- (50.5,71) node[midway,yshift=10pt] {$q$};  % 42.5→45.5，47.5→50.5
    
    \path[line width=0.30mm, draw=L] (38,65) -- (40,70);  % 35→38，37→40
    \path[line width=0.30mm, draw=L] (46,70) circle (0.50mm);  % 43→46
    \path[line width=0.30mm, draw=L, fill=F] (47,70) circle (0.20mm);  % 44→47
    \path[line width=0.30mm, draw=L, fill=F] (48,70) circle (0.20mm);  % 45→48
    \path[line width=0.30mm, draw=L, fill=F] (49,70) circle (0.20mm);  % 46→49
    \path[line width=0.30mm, draw=L] (50,70) circle (0.50mm);  % 47→50
    \path[line width=0.30mm, draw=L] (48,65) -- (46,70);  % 45→48，43→46
    \path[line width=0.30mm, draw=L] (48,65) -- (50,70);  % 45→48，47→50
    \path[line width=0.30mm, draw=L] (66,70) circle (0.50mm);  % 63→66
    \path[line width=0.30mm, draw=L] (70,70) circle (0.50mm);  % 67→70
    \path[line width=0.30mm, draw=L, fill=F] (67,70) circle (0.20mm);  % 64→67
    \path[line width=0.30mm, draw=L, fill=F] (68,70) circle (0.20mm);  % 65→68
    \path[line width=0.30mm, draw=L, fill=F] (69,70) circle (0.20mm);  % 66→69
    \path[line width=0.30mm, draw=L] (76,70) circle (0.50mm);  % 73→76
    \path[line width=0.30mm, draw=L] (80,70) circle (0.50mm);  % 77→80
    \path[line width=0.30mm, draw=L, fill=F] (77,70) circle (0.20mm);  % 74→77
    \path[line width=0.30mm, draw=L, fill=F] (78,70) circle (0.20mm);  % 75→78
    \path[line width=0.30mm, draw=L, fill=F] (79,70) circle (0.20mm);  % 76→79
    
    \draw[decorate,decoration={brace,amplitude=5pt}] (65.5,71) -- (70.5,71) node[midway,yshift=10pt] {$q$};  % 62.5→65.5，67.5→70.5
    \draw[decorate,decoration={brace,amplitude=5pt}] (75.5,71) -- (80.5,71) node[midway,yshift=10pt] {$p$};  % 72.5→75.5，77.5→80.5
    
    \path[line width=0.30mm, draw=L] (68,65) -- (66,70);  % 65→68，63→66
    \path[line width=0.30mm, draw=L] (68,65) -- (70,70);  % 65→68，67→70
    \path[line width=0.30mm, draw=L] (78,65) -- (76,70);  % 75→78，73→76
    \path[line width=0.30mm, draw=L] (78,65) -- (80,70);  % 75→78，77→80
    \path[line width=0.30mm, draw=L] (38,65) -- (48,65);  % 35→38，45→48
    \path[line width=0.30mm, draw=L] (38,65) -- (28,80);  % 35→38，25→28
    \path[line width=0.30mm, draw=L] (28,80) -- (23,82);  % 25→28，20→23
    \path[line width=0.30mm, draw=L] (28,80) -- (23,78);  % 25→28，20→23
    \path[line width=0.30mm, draw=L] (38,65) -- (28,50);  % 35→38，25→28
    \path[line width=0.30mm, draw=L] (28,50) -- (23,52);  % 25→28，20→23
    \path[line width=0.30mm, draw=L] (28,50) -- (23,48);  % 25→28，20→23
\end{tikzpicture}
    \caption{a star-like path tree $SLP(p,q;c;d,e)$}
    \label{fig:slp_structure}
\end{figure}
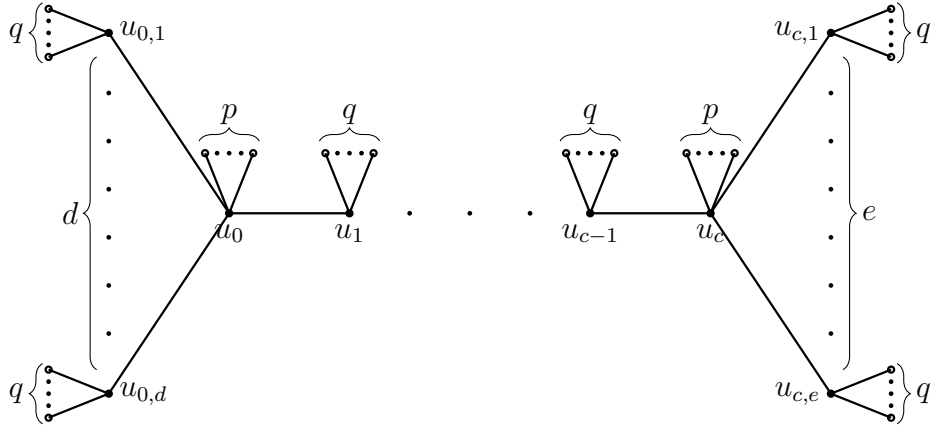

\begin{definition}
    Let $q \geq 1$ and $p, c,d,e \geq 0$ be integers.
    The star-like path tree $SLP(p,q;c;d,e)$ is obtained through the following process. 
    Let $u_0 \sim u_1 \sim \cdots \sim u_c$ be a path of length $c$. 
    Attach $d$ vertices to $u_0$, which are denoted by $u_{0,i}$ ($1\leq i\leq d$), and attach $e$ vertices to $u_c$, which are denoted by $u_{c,j}$ ($1\leq j\leq e$).
    Then attach $p$ leaves to $v \in \set{u_0, u_c}$, and attach $q$ leaves to each $u \in \{u_{0,i}|1\leq i\leq d\}\cup \{u_{c,j}|1\leq j\leq e\}\cup \{u_k|1\leq k\leq c-1\}$. 
    See~\cref{fig:slp_structure}.
\end{definition}

\begin{example}\label{exa:SLP-a-1}
    Let $a \geq 1$, $k \geq 2$, and $s \in [1,k-1]$ be integers. 
    Let $b = ak + s$. 
    Then $T = SLP(a + s, a; 0; k - 1, 0)$ is a tree with $b$ leaves as boundary and $k$ interior vertices. 
    The first Dirichlet eigenvalue of $T$ is given by $\lambda_1(T) = \sigma_1(a,k,s)$, where
    \begin{align}
        \sigma_1 = \sigma_1(a,k,s) = \frac{k + 2a + s - \sqrt{(k + s - 2)^2 + 4k - 4}}{2}.   
    \end{align}
    In particular, $a < \sigma_1 < a + 1$.
    The $\lambda_1$-eigenfunction is given by
    \begin{align*}
        f_1(u) = 
        \begin{cases} 
            1, &  u = u_0, \\
            \frac{1}{a + 1 - \sigma_1}, & \textotherwise.
        \end{cases}
    \end{align*}
\end{example}

\begin{example}\label{exa:SLP-(a+1)-(k-1)}
    Let $a \geq 0$ and $k \geq 2$ be integers. 
    Let $b = ak + k - 1$. 
    Then $T = SLP(a, a+1; 0; k - 1, 0)$ is a tree with $b$ leaves as boundary and $k$ interior vertices. 
    The first Dirichlet eigenvalue of $T$ is given by $\lambda_1(T) = \sigma_2(a,k)$, where
    \begin{align}
        \sigma_2 = \sigma_2(a, k) = \frac{k + 2a + 1 - \sqrt{(k - 1)^2 + 4}}{2}.  
    \end{align}
    In particular, $a < \sigma_2 < a + 1$.
    The $\lambda_1$-eigenfunction is given by
    \begin{align*}
        f_2(u) = 
        \begin{cases} 
            1, &  u = u_0, \\
            \frac{1}{a + 2 - \sigma_2}, & \textotherwise.
        \end{cases}
    \end{align*}
\end{example}

\subsection{\texorpdfstring{The trees with $b=ak+1$ or $b=ak+k-1$}{The trees with b=ak+1 or b=ak+k-1}}\label{subsec:b=ak+1}

Let $a \geq 1$ and $k \geq 2$ be integers. 
Next, we determine the extremal tree with maximum first Dirichlet eigenvalue among all trees with $b=ak+1$, or $b=ak+k-1$ leaves and $k$ interior vertices.

\begin{theorem}\label{the: SLP-a-1}
    Let $a \geq 1$ and $k \geq 2$ be integers and let $b = ak + 1$.
    Let $T$ be a tree with $k$ interior vertices and $b$ leaves. 
    Then
    \begin{align*}
        \lambda_1(T) \leq \sigma_1(a, k, 1) = \frac{k + 2a + 1 - \sqrt{(k - 1)^2 + 4k - 4}}{2}.   
    \end{align*}
    The equality holds if and only if $T$ is the star-like path tree $SLP(a+1, a;0; k - 1, 0)$.
\end{theorem}

\begin{proof}
    By the pigeonhole principle, there exists an interior vertex $v_0$ with exactly $a_{v_0} = |N(v_0) \cap B| \geq a + 1$ adjacent leaves. 
    Denote the number of interior vertices adjacent to $v_0$ by $m_{v_0} = |N(v_0) \cap \Omega|$.
    Then $m_{v_0} \leq k - 1$.
    Let $T_0$ be the tree $SLP(a+1, a;0; k - 1, 0)$. 
    Recall the $\sigma_1$-eigenfunction $f_1$ of $T_0$ in~\cref{exa:SLP-a-1} is given by 
    \begin{align*}
       f_1(u) = 
       \begin{cases} 
           1 & u = u_0, \\
           \frac{1}{a+1 - \sigma_1} & \text{otherwise},
        \end{cases}
    \end{align*}
    where $\sigma_1 = \sigma_1(a, k, 1)$.  
    Note that $f_1(u) > f_1(u_0) > 0$ for $u \neq u_0$. 
    Define a test function on $T$ by
    \begin{align*}
        g_1(v) = 
        \begin{cases} 
            1 & v = v_0, \\
            \frac{1}{a+1 - \sigma_1} & \text{otherwise}.
        \end{cases}
    \end{align*}
    We have
    \begin{align*}
        \lambda_1(T) &\leq R_T(g_1) \\
            &= \frac{m_{v_0} (1 - t_1)^2 + (ak+1-a_0) (t_1-0)^2 + a_0 (1-0)^2}{1^2 + (k-1) t_1^2} \\
            &= \frac{m_{v_0}(1 - t_1)^2 + (ak + 1 - a_0)t_1^2 + (a_0 - 1 - a) + (a + 1)}{1 + (k-1) t_1^2} \\
            &\leq \frac{(k-1)(1 - t_1)^2 + (ak + 1 - a_0)t_1^2 + (a_0 - 1 - a) t_1^2 + (a+1)}{1 + (k-1) t_1^2} \\
            &= \frac{(k-1)(1 - t_1)^2 + a (k-1) t_1^2 + (a + 1)}{1 + (k-1) t_1^2} = R_{T_0}(f_1) = \sigma_1,
    \end{align*}
    where $t_1 = \frac{1}{a+1-\sigma_1} > 1$. 
    If the equality holds, then $m_{v_0} = k - 1$, $a_0 = a + 1$. 
    If $m_{v_0} = k - 1$, $a_0 = a + 1$, then $g_1$ is an $\lambda_1(T)$-eigenfunction of $T$, and hence the equality holds. 
    Therefore, the tree $SLP(a+1, a;0; k - 1, 0)$ is the only extremal graph.
\end{proof}

\begin{theorem}\label{the: SLP-(a+1)-(k-1)}
    Let $a \geq 0$ and $k \geq 2$ be integers and let $b = ak + k - 1$.
    Let $T$ be a tree with $k$ interior vertices and $b$ leaves. 
    Then
    \begin{align*}
        \lambda_1(T) \leq \sigma_2(a, k) = \frac{k + 2a +1 - \sqrt{(k - 1)^2+4}}{2}.  
    \end{align*}
    The equality holds if and only if $T$ is the star-like path tree $SLP(a, a+1; 0; k- 1, 0)$.
\end{theorem}

\begin{proof}
    By~\cref{exa:SLP-(a+1)-(k-1)}, we have $f_2(u_0) > f_2(u)$ for any $u\neq u_0$. 
    The proof is similar to that of~\cref{the: SLP-a-1}.
\end{proof}

\subsection{\texorpdfstring{The trees with $b=ak+2$}{The trees with b=ak+2}}\label{subsec:b=ak+2}

Let $a \geq 1$ and $k \geq 2$ be integers. 
Next, we determine the extremal tree with maximum first Dirichlet eigenvalue among all trees with $b=ak+2$ leaves and $k$ interior vertices.
We first consider the cases where $k$ is even.

\begin{example}\label{exa:SLP-a-2-even}
    Let $a \geq 1$ and $k \geq 4$ be integers, and let $b = ak + 2$.
    Suppose $k$ is even. 
    Then $T = SLP(a+1, a; 1; \frac{k - 2}{2}, \frac{k - 2}{2})$, $SLP(a+1, a; 2; \frac{k - 4}{2}, \frac{k - 2}{2})$, $SLP(a+1, a; 3; \frac{k - 4}{2}, \frac{k - 4}{2})$ are trees with $b$ leaves as boundary and $k$ interior vertices. 
    And $\lambda_1(T) = \sigma_3(a, k)$, where
    \begin{align}
        \sigma_3 = \sigma_3(a, k) = \frac{k + 4a + 2 - \sqrt{(k + 2)^2-16}}{4}.   
    \end{align}
    There are two facts as follows.
    \begin{enumerate}
        \item $a< \sigma_3 < a + 1$;
        \item $\sigma_3$ is an irrational number.
    \end{enumerate}
    One $\lambda_1$-eigenfunction of $SLP(a + 1, a; 1; \frac{k - 2}{2}, \frac{k - 2}{2})$ is given by
    \begin{align*}
        f_{3,\text{I}}(u) = 
        \begin{cases} 
            1, &  u \in \set{u_0, u_1}, \\
            \frac{1}{a + 1 - \sigma_3}, & \textotherwise.
        \end{cases}
    \end{align*}
    One $\lambda_1$-eigenfunction of $SLP(a + 1, a; 2; \frac{k - 4}{2}, \frac{k - 2}{2})$ is given by
    \begin{align*}
       f_{3,{\text{II}}}(u) = 
       \begin{cases} 
           1 & u = u_{0, i}, u_1, u_2, \ 1 \leq i \leq \frac{k}{2} - 2, \\
           a + 1 - \sigma_3 & u = u_0, \\
           \frac{1}{a + 1 - \sigma_3} & u = u_{2, i},\ 1 \leq i \leq \frac{k}{2} - 1.
        \end{cases}
    \end{align*}
    One $\lambda_1$-eigenfunction of $SLP(a + 1, a; 3; \frac{k - 4}{2}, \frac{k - 4}{2})$ is given by
    \begin{align*}
        f_{3,\text{III}}(u) = 
        \begin{cases} 
            1, &  u \in \set{u_0, u_3}, \\
            \frac{1}{a + 1 - \sigma_3}, & \textotherwise.
        \end{cases}
    \end{align*}
\end{example}

\begin{theorem}\label{the:SLP-a-2-even}
    Let $a \geq 1$ and $k \geq 4$ be integers, and let $b = ak + 2$. 
    Let $T$ be a tree with $k$ interior vertices and $b$ leaves. 
    Suppose $k$ is even, then
    \begin{align*}
        \lambda_1(T) \leq \sigma_3(a, k) = \frac{k + 4a + 2 - \sqrt{(k + 2)^2-16}}{4}.
    \end{align*}
    The equality holds if and only if $T$ is the star-like path tree $SLP(a+1,a;1;\frac{k-2}{2},\frac{k-2}{2})$, $SLP(a+1,a;2;\frac{k-4}{2},\frac{k-2}{2})$, or $SLP(a+1,a;3;\frac{k-4}{2},\frac{k-4}{2})$. 
\end{theorem}

\begin{proof}
    Let $T_1, T_2, T_3$ be the trees $SLP(a+1,a;1;\frac{k-2}{2},\frac{k-2}{2})$, $SLP(a+1,a;2;\frac{k-4}{2},\frac{k-2}{2})$, and $SLP(a+1,a;3;\frac{k-4}{2},\frac{k-4}{2})$ respectively.
    By~\cref{exa:SLP-a-2-even}, we have $\lambda_1(T_1) = \lambda_1(T_2) = \lambda_1(T_3) = \sigma_3(a, k)$.
    Suppose $\lambda_1(T) \geq \sigma_3(a, k)$ and $T$ is an extremal tree.
    For each interior vertex $v \in \Omega$, let $a_v = |N(v) \cap B|$ be the number of leaves adjacent to $v$.
    By the pigeonhole principle, there exist two interior vertices $v_0$ and $v_1$ such that $a_{v_0} + a_{v_1} \geq 2a + 2$.
    Let $m_{v_0} = |N(v_0) \cap (\Omega \setminus \set{v_0, v_1})|$ and $m_{v_1} = |N(v_1) \cap (\Omega \setminus \set{v_0, v_1})|$ denote the number of interior vertices, except $v_0$ and $v_1$, adjacent to $v_0$ and $v_1$ respectively. 
    Then $m_{v_0} + m_{v_1} \leq k - 1$, and the equality may hold if $v_0$ and $v_1$ share a common neighbor.
    We divide the discussion into two cases. 
    \begin{enumerate}[{label=(\arabic*)}]
        \item\label{SLP-a-1-even-case1} $N(v_0)\cap N(v_1) = \emptyset$.

        In this case, $m_{v_0} + m_{v_1} \leq k - 2$.
        Recall the $\lambda_1$-eigenfunction of $SLP(a+1,a;1;\frac{k-2}{2},\frac{k-2}{2})$ is given by 
        \begin{align*}
           f_{3,\text{I}}(u) &= 
           \begin{cases} 
               1, & v = u_0, u_1, \\
               \frac{1}{a+ 1 - \sigma_3}, & \text{otherwise},
           \end{cases}
        \end{align*}
        where $\sigma_3 = \sigma_3(a, k)$. 
        Note that $f_3(u) >  f_3(u_0)$ for $u \neq u_0, u_1$. 
        Define a test function on $T$ by
        \begin{align*}
            g_3(v) = 
            \begin{cases} 
               1 & v = v_0, v_1, \\
               \frac{1}{a + 1 - \sigma_3} & \text{otherwise}.
            \end{cases}
        \end{align*}
        We have
        \begin{align}
            \lambda_1(T) &\leq R_T(g_3) \nonumber \\
                &= \frac{(m_{v_0}+m_{v_1})(1 - t_3)^2 + (ak + 2 - a_{v_0} -a_{v_1})t_3^2 + (a_{v_0}+a_{v_1})}{2 + (k - 2)t_3^2} \nonumber \\
                &\leq \frac{(k - 2)(1 - t_3)^2 + [(ak + 2 - a_{v_0} - a_{v_1}) + (a_{v_0} + a_{v_1} - 2a - 2)] t_3^2 + 2(a + 1)}{2 + (k - 2) t_3^2} \nonumber \\
                &= \frac{(k - 2)(1 - t_3)^2 + a(k - 2)t_3^2 + 2(a + 1)}{2 + (k - 2)t_3^2} \nonumber \\
                &= R_{T_1}(f_{3, \text{I}}) = \sigma_3, \label{eq:case1}
        \end{align}
        where $t_3 = \frac{1}{a+ 1 - \sigma_3} > 1$. 
        Note that the equality holds if $T$ is  $SLP(a+1,a;1;\frac{k-2}{2},\frac{k-2}{2})$ or $SLP(a+1,a;3;\frac{k-4}{2},\frac{k-4}{2})$. 
        If the equality holds, then $m_{v_0} + m_{v_1} = k - 2$, $a_{v_0} + a_{v_1} = 2a+2$, and $g_3$ is a $\lambda_1$-eigenfunction of $T$. 
        Suppose that $T$ is an extremal tree. 
        \begin{enumerate}
            \item\label{d(v_0 v_1) = 1} Suppose $d(v_0, v_1) = 1$.  
            
            Let $w \in (N(v_0) \cup N(v_1)) \setminus \set{v_0, v_1}$ be an interior vertex adjacent to one of $v_0$ and $v_1$. 
            The function $g_3$ is a $\sigma_3$-eigenfunction implies that 
            \begin{align*}
                \sigma_3 g_3(w) = \sum \limits_{y \sim w}(g_3(w)-g_3(y)).
            \end{align*}
            Then
            \begin{align}\label{equ:d_w=a+1}
                \frac{1}{d_{w} - \sigma_3} = \frac{1}{a + 1 - \sigma_3},
            \end{align}
            where $d_w-\sigma_3 \neq 0$ since $\sigma_3$ is irrational.
            Hence, we have $d_{w} = a+1$.
            Since $g_3$ is an eigenfunction and $g_3(v_0) = g_3(v_1)$, we have $(a + 1 - \sigma_3)(m_{v_0} + a_{v_0} - m_{v_1} - a_{v_1}) = m_{v_0} - m_{v_1}$. 
            Note that $\sigma_3$ is an irrational number when $k \geq 4$.
            Therefore, $m_{v_0}-m_{v_1}=m_{v_0}+a_{v_0}-m_{v_1}-a_{v_1}=0$, which implies that $m_{v_0} = m_{v_1} = (k-2)/2$ and $a_{v_0} = a_{v_1} = a+1$.
            Therefore, the unique extremal tree in this subcase is $SLP(a+1,a;1;\frac{k-2}{2},\frac{k-2}{2})$.
            
            \item Suppose $d(v_0, v_1)\geq 3$.
            
            Since $m_{v_0} + m_{v_1} = k - 2$, we have $d(v_0,v_1)=3$.  
            Denote the unique path connecting $v_0$ and $v_1$ by $v_0 \sim x_0 \sim x_1 \sim v_1$. 
            Let $w = (N(v_0) \cup N(v_1)) \setminus \set{v_0, v_1, x_0, x_1}$ be an interior vertex adjacent to either $v_0$ or $v_1$, which is different from $x_0$ and $x_1$. 
            A similar argument as above tells us $d_{w} = a + 1$, $m_{v_0} = m_{v_1} = \frac{k-2}{2}$, and $a_{v_0}=a_{v_1}=a + 1$. 
            Since $g_3$ is an eigenfunction, we have 
            \begin{align*}
                \sigma_3 g_3(x_i)=\sum \limits_{y \sim x_i}(g_3(x_i)-g_3(y)).
            \end{align*}
             Then
            \begin{align*}
               \frac{1}{d_{x_i}-\sigma_3-1}=\frac{1}{a+ 1 - \sigma_3},
            \end{align*}
            where $d_{x_i}-\sigma_3-1 \neq 0$ since $\sigma_3$ is irrational.
            Then $d_{x_0}=d_{x_1}=a+2$.
            Therefore, the unique extremal graph in this subcase is $SLP(a+1,a;3;\frac{k-4}{2},\frac{k-4}{2})$.
        \end{enumerate}
        
        \item $N(v_0)\cap N(v_1) \neq \emptyset$.
        
        In this case we have $d(v_0,v_1)=2$.
        Denote the unique path connecting $v_0$ and $v_1$ by $v_0 \sim x \sim v_1$. 
        
        \begin{claim*}
            Every interior vertex $v \in \Omega \setminus \set{v_0, v_1}$ is adjacent to at least one of $v_0$ and $v_1$. 
            And hence $m_{v_0}+m_{v_1}=k-1$. 
        \end{claim*}
        
        \begin{proof}
            Suppose there exists at least one interior vertex $v \in \Omega \setminus \set{v_0, v_1}$ that is not adjacent to $v_0$ or $v_1$. 
            Then $m_{v_0} + m_{v_1} \leq k - 2$, and~\ref{eq:case1} still holds.
            If the equality holds, then $m_{v_0} + m_{v_1} = k - 2$, $a_{v_0} + a_{v_1} = 2 a + 2$, and $g_3$ is the eigenfunction of $\lambda_1(T)$. 
            Hence, there is exactly one interior vertex $z$ is not adjacent to either of $v_0$ and $v_1$.
            Since $g_3$ is an eigenfunction, we have 
            \begin{align*}
                \sigma_3 g_3(z)=\sum \limits_{y \sim z}(g_3(z)-g_3(y)).
            \end{align*}
            And we further deduce that $\sigma_3 = d_{z}-1$.
            This contradicts the fact that $\sigma_3$ is an irrational number. 
            The claims follows from $N(v_0) \cap N(v_1) = \set{x}$ and $N(v_0) \cup N(v_1) = \Omega \setminus \set{v_0, v_1}$. 
        \end{proof}
        
        \begin{claim*}
            Each interior vertex $v \in \Omega \setminus \set{v_0, v_1}$ is adjacent to exactly $a$ leaves in the extremal tree.
            Namely, $a_v = a$ for $v \in \Omega \setminus \set{v_0, v_1}$.
        \end{claim*}
        
        \begin{proof}
            If there exists $v \in \Omega \setminus \set{v_0, v_1}$ such that $a_v \leq a-1$. 
            By~\cref{lambda at most minimum degree}, we have $\lambda_1(T)\leq a < \sigma_3(a, k)$, and $T$ can not be extremal.
            Therefore, for each interior vertex $v \in \Omega \setminus \set{v_0, v_1}$, we have $a_v \geq a$. 
            Since $\sum_{v \in \Omega} a_v = ak+2$ and $a_{v_0} + a_{v_1} \geq 2a+2$, we conclude that $a_v = a$ for $a \in \Omega \setminus \set{v_0, v_1}$.
        \end{proof}
        
        \begin{claim*}
            Both of $v_0$ and $v_1$ are adjacent to exactly $a + 1$ leaves. 
            Namely, $a_{v} = a+1$ for $v \in \set{v_0, v_1}$. 
        \end{claim*}
        
        \begin{proof}
            Note that $a_{v_0} + a_{v_1} \geq 2a+2$. 
            Without loss of generality, assume that $a_{v_0} \geq a+2$. 
            Consider $v' \in \Omega \setminus \set{v_0, v_1}$. 
            By the above claim we have $a_{v'} = a$.
            Hence, $a_{v_0} + a_{v'} \geq 2a+2$ and $N(v_0) \cap N(v') = \emptyset$, which reduces to~\ref{SLP-a-1-even-case1}. 
            Therefore, $a_{v_0} = a_{v_1} = a+1$. 
        \end{proof}
        
        Next, we need to determine $m_{v_0}$ and $m_{v_1}$.
        Without loss of generality, we assume that $m_{v_1} \geq m_{v_0}$.
        Recall one $\sigma_3$-eigenfunction of $T_2$ isgiven by
        \begin{align*}
            f_{3,{\text{II}}}(u) = 
            \begin{cases} 
                1 & u = u_{0, i}, u_1, u_2, \ 1 \leq i \leq \frac{k}{2} - 2, \\
                a+ 1 - \sigma_3 & u = u_0, \\
                \frac{1}{a+ 1 - \sigma_3} & u = u_{2, i},\ 1 \leq i \leq \frac{k}{2} - 1,
            \end{cases}
        \end{align*}
        where $\sigma_3=\sigma_3(a, k)$, which is shown in~\cref{exa:SLP-a-2-even}. 
        It is clear that $f_{3,{\text{II}}}(u_{0, i}) =  f_{3,{\text{II}}}(u_1) = f_{3,{\text{II}}}(u_2) >f_{3,{\text{II}}}(u_0)$ and $(f_{3,{\text{II}}}(u_{0, i}) - f_{3,{\text{II}}}(u_0))^2 > (f_{3,{\text{II}}}(u_{0, i}) - f_{3,{\text{II}}}(u_{2}))^2 = 0 $. 
        
        Let $v_{0, i} \in N(v_0) \cap (\Omega \setminus \set{x}), 1 \leq i \leq m_{v_0} - 1$ be the neighbors of $v_0$ other than $x$, and let $v_{1i} \in N(v_1) \cap (\Omega \setminus \set{x}), 1 \leq i \leq m_{v_1} - 1$ be the neighbors of $v_1$ other than $x$.
        Define a test function on $T$ by
        \begin{align*}
            g_{3,{\text{II}}}(v) = 
            \begin{cases} 
                1 & v = v_{0,i}, x, v_1, v_{1,j}, \ 1 \leq i \leq m_{v_0}-1,
                \frac{k}{2} \leq j \leq m_{v_1}-1,\\
                a+ 1 - \sigma_3 & v = v_0, \\
                \frac{1}{a+ 1 - \sigma_3} & v = v_{1,i},\ 1 \leq i \leq \frac{k}{2} - 1.
            \end{cases}
        \end{align*}
        Then 
        \begin{align*}
            \lambda_1(T) \leq& R_T(g_{3, \text{II}}) \\
            =& \frac{(m_{v_0}-1) (t_3^{-1} - 1)^2 + (t_3^{-1} - 1)^2 +(1 - 1)^2}{\frac{k}{2} 1^2 + t_3^{-2} + (\frac{k}{2}-1) t_3^2} \\
            &+ \frac{\frac{k-2}{2} (1 - t_3)^2 + (m_{v_1}-\frac{k-2}{2} -1) (1 - 1)^2}{\frac{k}{2} 1^2 + t_3^{-2} + (\frac{k}{2}-1) t_3^2} \\
            &+ \frac{(a+1) t_3^{-2} + \frac{ka+2}{2} 1^2 +\frac{a(k-2)}{2} t_3^2}
            {\frac{k}{2} 1^2 + t_3^{-2} + (\frac{k}{2}-1) t_3^2} \\
            \leq& \frac{\frac{k - 4}{2} (t_3^{-1} - 1)^2 + (t_3^{-1} - 1)^2 + (1 - 1)^2}{\frac{k}{2} 1^2 + t_3^{-2} + (\frac{k}{2} - 1) t_3^2} \\
            &+ \frac{\frac{k - 2}{2} (1 - t_3)^2 + (a+1) t_3^{-2}+\frac{ka+2}{2} 1^2+\frac{a(k-2)}{2}t_3^2}
            {\frac{k}{2} 1^2 + t_3^{-2} + (\frac{k}{2} - 1) t_3^2} \\
            =& R_{T_2}(f_{3,\text{II}}) = \sigma_3, 
        \end{align*}
        where $t_3 = \frac{1}{a+ 1 - \sigma_3}$.
        The equality holds if and only if $m_{v_0} - 1 = \frac{k-4}{2}$, $m_{v_1} - 1 = \frac{k-2}{2}$, and $g_{3,{\text{II}}}$ is an eigenfunction. 
        Therefore, the unqiue extremal tree in this case is  $SLP(a+1,a;2;\frac{k-4}{2},\frac{k-2}{2})$.
    \end{enumerate}
    The proof is complete.
\end{proof}

Next we consider the cases where $k$ is odd.

\begin{example}\label{exa:SLP-r-2-odd}
    For integers $a \geq 1$ and $k \geq 5$, let $b= ak +2$. 
    If $k$ is odd, then $SLP(a + 1,  a; 2; \frac{k-3}{2}, \frac{k-3}{2})$ is a tree with $b$ leaves as boundary and $k$ interior vertices. 
    The Dirichlet eigenvalues of the tree $SLP(a + 1,  a; 2; \frac{k-3}{2}, \frac{k-3}{2})$ are roots of $F(\lambda) = (\lambda - (a + 1))^{k-5} Q(\lambda) R(\lambda)$, where
    \begin{align*}
       Q(\lambda) &= \lambda^3 - \lambda^2( \frac{k}{2} + 3a + \frac{7}{2}) + \lambda ( ka + k + 3a^2 + 7a + 3)\\
       &-( \frac{ka^2}{2} + ka + a^3 + \frac{7a^2}{2} + 3a + 2),
    \end{align*}
    and 
    \begin{align*}
       R(\lambda) = 2\lambda^2 - (k + 4a + 3)\lambda + (ka + 2a^2 + 3a + 4).
    \end{align*}
    %First, we consider the roots of $Q(\lambda)$. 
    %Here $Q(a) = Q(a+2) = -2$, $Q(a + 1) = \frac{1}{2}k - \frac{3}{2}$ and $Q(+\infty)>0$.
    %Thus, we find the distribution of three roots of $Q(\lambda)$, which implies that these roots are larger than $a$.
    %Second, we consider the roots of $R(\lambda)$.
    %For the equation $R(\lambda) = 0$, the roots are  
    %$\frac{k+4a+3-\sqrt{k^2 + 6k - 23}}{4}$ and $\frac{k+4a+3+\sqrt{k^2 + 6k - 23}}{4}$,  
    %where $\frac{k+4a+3-\sqrt{k^2 + 6k - 23}}{4} > a$.
    %Then $Q(\frac{k+4a+3-\sqrt{k^2 + 6k - 23}}{4}) = \frac{\sqrt{k^2 + 6k - 23} + 1 - k}{2}$.
    %Therefore, $\lambda_1(SLP(a + 1, a; 2; \frac{k-3}{2}, \frac{k-3}{2})) > a$.
    One $\lambda_1$-eigenfunction of $SLP(a + 1, a; 2; \frac{k-3}{2}, \frac{k-3}{2})$ is given by
   \begin{align*}
       f_4(u) = 
       \begin{cases} 
           1 & u = u_0, u_2, \\ 
           \frac{2}{a + 2 - \sigma_4} & u = u_1, \\
           \frac{1}{a + 1 - \sigma_4} & \text{otherwise}, 
        \end{cases}
    \end{align*}
    where $a < \sigma_4(a, k) = \lambda_1(SLP(a + 1, a; 2; \frac{k-3}{2}, \frac{k-3}{2})) < a+1$ and $\sigma_4(a, k)$ is a root of $Q(\lambda)$.
\end{example}

%\begin{example}\label{exa:SLP-(a + 1)-(k-2)-odd}
    %For integers $a \geq 0$, $ k \geq 5$, let $b=ak+ k-2$.
    %If $k$ is odd, then $SLP(a, a + 1; 2; \frac{k-3}{2}, \frac{k-3}{2})$ is a tree with $b$ leaves as boundary and $k$ interior vertices. 
    %The Dirichlet eigenvalues of the tree $SLP(a, a + 1; 2; \frac{k-3}{2}, \frac{k-3}{2})$ are roots of $F(\lambda) = (\lambda - (a+2))^{k-5} \cdot Q(\lambda) \cdot R(\lambda)$, where
   %\begin{align*}
        %Q(\lambda) &= 
        %\lambda^3 + \lambda^2(-\frac{k}{2} - 3a - \frac{9}{2})\\
        %&+ \lambda(ka + 2k + 3a^2 + 9a + 3) + \frac{-ka^2 - 4ka - 3k - 2a^3 - 9a^2 - 6a + 5}{2},
   %\end{align*}
   %and 
   %\begin{align*}
       %R(\lambda) = \lambda^2 + (-\frac{k}{2} - 2a - \frac{3}{2})\lambda + %\frac{ka + k + 2a^2 + 3a + 1}{2}.
   %\end{align*}
   %One $\lambda_1$-eigenfunction of $SLP(a, a + 1; 2; \frac{k-3}{2}, \frac{k-3}{2})$ is given by
   %\begin{align*}
       %f_6(u) = 
       %\begin{cases} 
           %1 & u = u_0, u_2, \\ 
           %\frac{2}{a + 3 - \sigma_6} & u = u_1, \\
           %\frac{1}{a + 2 - \sigma_6} & \text{otherwise}, 
        %\end{cases}
    %\end{align*}
    %\textcolor{red}{where $a < \sigma_6 = \sigma_6(a, k) = \lambda_1(SLP(a, a + 1; 2; \frac{k-3}{2}, \frac{k-3}{2})) < a+1$ and $\sigma_6(a, k)$ is a root of $Q(\lambda)$.}
%\end{example}

\begin{theorem}\label{thm:SLP-a-2-odd}
    Let $a \geq 1$ and $k \geq 5$ be integers, and let $b = ak + 2$. 
    Let $T$ be a tree with $k$ interior vertices and $b$ leaves. 
    Suppose $k$ is odd, then
    \begin{align*}
        \lambda_1(T) \leq \lambda_1(SLP(a + 1,a; 2;\frac{k-3}{2},\frac{k-3}{2})).
    \end{align*}
    The equality holds if and only if $T$ is the star-like path tree $SLP(a + 1,a;2;\frac{k-3}{2},\frac{k-3}{2})$.
\end{theorem}

\begin{proof}
    Let $T_0$ be the tree $SLP(a + 1,a;2;\frac{k-3}{2},\frac{k-3}{2})$. 
    By~\cref{exa:SLP-r-2-odd}, we have $\lambda_1(T_0) = \sigma_4(a, k)$. 
    Suppose that $T$ is an extremal tree. 
    Then $\lambda_1(T) \geq \sigma_4(a, k) > a$. 
    For each interior vertex $v \in \Omega$, let $a_v = |N(v) \cap B|$ be the number of leaves adjacent to $v$.
    
    \begin{claim*}
         There exist three distinct interior vertices $v_0, v_1, x \in \Omega$ such that $a_{v_0} + a_{v_1} \geq 2a + 2$, $a_{v_0} \geq a+1$, $a_x \geq a$.
    \end{claim*}
    
    \begin{proof}
        By the pigeonhole principle, there exist two interior vertices $x_0$ and $x_1$ such that $a_{x_0} + a_{x_1} \geq 2a + 2$.
        Without loss of generality, suppose that $a_{x_0} \geq a_{x_1}$.
        Then $a_{x_0} \geq a+1$.
        
        If $d(x_0, x_1) = k-1$, then all interior vertices are on a path of length $k-1$ with two ends $x_0$ and $x_1$. 
        If there exists $x_2 \in \Omega \setminus \set{x_0, x_1}$ such that $a_{x_2} \geq a$, then the claim follows by taking $v_0 = x_0, v_1 = x_1$, and $x = x_2$. 
        Otherwise, we have $a_{v} \leq a-1$ for every $v \in \Omega \setminus \set{x_0, x_1}$.
        By~\cref{lambda at most minimum degree}, we have $a_{x_1} \geq a$ and $a_v \geq a-1$ for every $v \in \Omega \setminus \set{x_0, x_1}$. 
        Therefore, $a_{v} = a-1$ for every $v \in \Omega \setminus \set{x_0, x_1}$.
        Since $\sum_{v \in \Omega} a_v = ak+2$, $a_{v} = a-1$ for every $v \in \Omega \setminus \set{x_0, x_1}$, and $a_{x_0} \geq a_{x_1}$, we have $a_{x_0} + a_{v} \geq \frac{ak+2 - (k-2)(a-1)}{2} +a -1 = a + \frac{k}{2} +a -1$ for any $v \in \Omega \setminus \set{x_0, x_1}$.
        The claim follows by taking $v_0 = x_0, v_1 = v'$, and $x = x_1$ for some $v' \in \Omega \setminus \set{x_0, x_1}$. 
        
        If $d(x_0, x_1) \neq k-1$, then there exists $x_2 \in \Omega \setminus \set{x_0, x_1}$ such that $|N(x_2) \cap \Omega| = 1$. 
        By~\cref{lambda at most minimum degree}, we have $a_{x_2} \geq a$. 
        The claim follows by taking $v_0 = x_0, v_1 = x_1$, and $x = x_2$.
    \end{proof}
    
    For each interior vertex $v \in \Omega$, let $m_{v} = |N(v) \cap (\Omega \setminus \set{v_0, v_1, x})|$. 
    Please keep in mind that $m_v$ depends on $v_0, v_1$, and $x$ implicitly. 
    Then $m_{v_0} + m_{v_1} \leq k - 2$, and the equality may hold if $v_0$ and $v_1$ share a common neighbor which is different from $x$.
    Let $m_{x}' = |N(x) \cap \set{v_0, v_1}| \leq 2$. 
    We divide the discussion into two cases.
    
    \begin{enumerate}[{label = (\arabic*)}]
        \item\label{itm:leq_k-3} $m_{v_0} + m_{v_1} \leq k-3$. 
            
            Recall that a $\lambda_1$-eigenfunction of $SLP(a + 1,a;2;\frac{k-3}{2},\frac{k-3}{2})$ is given by 
            \begin{align*}
                f_4(u) = 
                \begin{cases} 
                    1, & v = u_0, u_2, \\
                    s_4, & v = u_1, \\
                    t_4, & \text{otherwise},
                \end{cases}
            \end{align*}
            where $\sigma_4 = \sigma_4(a, k)$, $t_4 = \frac{1}{a + 1 - \sigma_4}$, and $s_4 = \frac{2}{a+2-\sigma_4}$. 
            Note that $1 < s_4 < t_4$ and $(t_4 - s_4)^2 < (s_4 - 1)^2 < (t_4 - 1)^2$ for every $u \in \Omega \setminus \set{u_0, u_1, u_2}$. 
            Define a test function on $T$ by
            \begin{align*}
                g_4(v) = 
                \begin{cases} 
                    1, & v = v_0, v_1, \\
                    s_4, & v = x, \\
                    t_4, & \text{otherwise}.
                \end{cases}
            \end{align*}
            
            \begin{claim*}
                We have $\lambda_1(T) \leq R_T(g_5) \leq \sigma_4(a, k)$. 
                If the equality holds, then $m_{v_0} + m_{v_1} = k - 3$, $m_x=0$, $m_x'=2$, $a_{v_0}+a_{v_1} = 2a+2$, $a_{x}=a$, and $g_4$ is a $\lambda_1$-eigenfunction of $T$.
            \end{claim*}
            
            \begin{proof}
                Since $a_x \geq a$, we have
                \begin{align*}
                    &(a_{v_0} + a_{v_1}) + a_x s_4^2 + [ak + 2 - (a_{v_0} + a_{v_1}) - a_x] t_4^2 \\
                    =&[a_{v_0} + a_{v_1} - 2(a + 1)] + 2(a + 1) + (a_x - a) s_4^2 + a s_4^2 \\
                    &+ [ak + 2 - (a_{v_0} + a_{v_1}) - a_x] t_4^2 \\
                    \leq& [a_{v_0} + a_{v_1} - 2(a + 1)] t_4^2 + 2(a + 1) + (a_x - a) t_4^2 + a s_4^2 \\
                    &+ [ak + 2 - (a_{v_0} + a_{v_1}) - a_x] t_4^2\\
                    =& 2(a + 1) + a s_4^2 + (ak - 3a) t_4^2.
                \end{align*}
                Since $m_{v_0} + m_{v_1} \leq k - 3$, we have
                \begin{align*}
                    &(m_{v_0}+m_{v_1})(t_4 - 1)^2 + m_x (t_4 - s_4)^2 + m_x' (1 - s_4)^2 \\
                    =& [m_{v_0} + m_{v_1} - (k - 3)] (t_4 - 1)^2 + (k - 3)(t_4 - 1)^2 + m_x (t_4 - s_4)^2 \\
                    &+ (m_x' - 2)(1 - s_4)^2 + 2(1 - s_4)^2 \\
                    \leq& (k - 3)(t_4 - 1)^2 + 2(1 - s_4)^2 + [m_{v_0} + m_{v_1} - (k - 3)] (t_4 - s_4)^2 \\
                    &+ m_x (t_4 - s_4)^2 + (m_x' - 2)(t_4 - s_4)^2 \\
                    \leq& (k - 3)(t_4 - 1)^2 + 2(1 - s_4)^2. 
                \end{align*}
                The last inequality holds because $m_{v_0} + m_{v_1} +m_x+m_x' \leq k-1$.
                Therefore,
                \begin{align}
                    \lambda_1(T) &\leq R_T(g_4) = \frac{(a_{v_0}+a_{v_1})  + a_{x} s_4^2 + (ak + 2 - (a_{v_0}+a_{v_1}) - a_{x}) t_4^2}{2  + s_4^2 + (k - 3) t_4^2} \nonumber \\
                    &\quad + \frac{(m_{v_0}+m_{v_1})(t_4 - 1)^2 + m_x (t_4 - s_4)^2 + m_x' (1 - s_4)^2}{2  + s_4^2 + (k - 3) t_4^2} \nonumber \\
                    &\leq \frac{2(a + 1) + a f_4^2(u_1) + (ak - 3a) t_4^2 +(k - 3) (t_4 - 1)^2 + 2 (f_4(u_1) - 1)^2}{2  + f_4^2(u_1) + (k - 3) t_4^2} \nonumber \\
                    &= R_{T_0}(f_4) = \sigma_4(a, k).  \label{eq:SLP(a-2-odd)}
                \end{align}
                If the equality holds, then $m_{v_0} + m_{v_1} = k - 3$, $m_x=0$, $m_x'=2$, $a_{v_0}+a_{v_1} = 2a+2$, $a_{x}=a$, and $g_4$ is a $\lambda_1$-eigenfunction of $T$.
            \end{proof}
            
            Since $m_x' = 2$, we have that $d(v_0, v_1) = 2$ and the path connecting $v_0$ and $v_1$ is $v_0 \sim x \sim v_1$.
            Since $m_{v_0} + m_{v_1} = k-3$ and $m_x = 0$, we have that every $v \in \Omega \setminus \set{v_0, x, v_1}$ is adjacent to $v_0$ or $v_1$. 
            Note that $a_x = a$.
            
            \begin{claim*}
                For every $v \in \Omega \setminus \set{v_0, x, v_1}$, we have $a_v = a$. 
            \end{claim*}
            
            \begin{proof}
                The function $g_4$ is a $\lambda_1$-eigenfunction of $T$ implies that 
                \begin{align*}
                    \sigma_4 g_4(v) = \sum \limits_{w: v \sim w}(g_4(v)-g_4(w)).
                \end{align*}
                Then
                \begin{align}\label{equ:odd-d_w=a+1}
                    \frac{1}{a_{v} + 1 - \sigma_4} = \frac{1}{a + 1 - \sigma_4}.
                \end{align}
                Since $d_v - \sigma_4 \neq 0$ we have $a_{v} = a$.
            \end{proof}
            
            \begin{claim*}\label{cla:SLP-2-odd-claim 2}
                We have $a_{v_0} = a_{v_1} = a + 1$.
            \end{claim*}
            
            \begin{proof}
                Note that $a_{v_0} + a_{v_1} = ak+2 - (k-2) a = 2a+2$. 
                We only need to show that $a_{v_0} = a+1$. 
                Suppose $a_{v_0} \geq a+2$. 
                Take $v' \in \Omega \setminus \set{v_0, x, v_1}$. 
                Then $a_{v'} = a$ and $a_{v_0} + a_{v'} \geq 2a+2$. 
                Take $v_0' = v_0$, $v_1'= v'$, and $x'=x$. 
                Then $m_{v_0'} + m_{v_1'} \leq k-3$. 
                Here $m_{v'} = |N(v') \cap (\Omega \setminus \set{v_0, v_1', x'})|$. 
                The same convention applies in the rest of the proof. 
                Then~\ref{eq:SLP(a-2-odd)} still holds. 
                But $d(v_0', v_1') \neq 2$, contradiction. 
                Therefore, $a_{v_0} = a_{v_1} = a + 1$.
            \end{proof}
            
            \begin{claim*}
                We have $m_{v_0} = m_{v_1} = \frac{k-3}{2}$. 
            \end{claim*}
            
            \begin{proof}
                Since $g_4$ is $\lambda_1$-eigenfunction, we have
                \begin{align*}
                    \sigma_4 g_4(v_i)=\sum \limits_{y: v_i \sim y}(g_4(v_i)-g_4(y))
                \end{align*}
                for $i = 0, 1$.
                Then
                \begin{align*}
                    m_{v_0} + 1 - (m_{v_0} t_4 + s_4) = m_{v_1} + 1 - (m_{v_1} t_4 + s_4).
                \end{align*}
                It reduces to $m_{v_0}(1 - t_4) = m_{v_1}(1 - t_4)$.
                Therefore, $m_{v_0} = m_{v_1} = \frac{k-3}{2}$.
            \end{proof}
            
            Combining the above three claims and we conclude that $T$ is the star-like path tree $SLP(a + 1,a; 2;\frac{k-3}{2},\frac{k-3}{2})$.
        \item $m_{v_0} + m_{v_1} = k-2$.
            
            In this case we have $d(v_0, v_1) = 2$. 
            Let the path connecting $v_0$ and $v_1$ be $v_0 \sim y \sim v_1$. 
            Since $m_{v_0} + m_{v_1} = k-2$, every $v \in \Omega \setminus \set{v_0, v_1, x}$ is adjacent to $v_0$ or $v_1$ or both. 
            
            If $x$ is adjacent to $y$, then
            every $v \in \Omega \setminus \set{v_0, y, v_1, x}$ is  adjacent to $v_0$ or $v_1$. 
            By~\cref{lambda at most minimum degree}, we have $a_{v} \geq a$.
            Let $v_0' = v_0, v_1' = v_1$, and $x' = v'$ for some $v' \in \Omega \setminus \set{v_0, y, v_1, x}$. 
            Then $m_{v_0'} + m_{v_1'} \leq k-3$, and it reduces to~\ref{itm:leq_k-3}.
            
            If $x$ is adjacent to $v_0$ or $v_1$, then every $v \in \Omega \setminus \set{v_0, v_1}$ is adjacent to $v_0$ or $v_1$ or both. 
            By~\cref{lambda at most minimum degree}, we have $a_v \geq a$ for $v \in \Omega \setminus \set{v_0, y, v_1}$.
            Then $a_{v_0} + a_y + a_{v_1} \leq ak+2 - (k-3)a = 3a+2$. 
            If $a_y \geq a$, then we take $v_0' = v_0, v_1' = v_1$, and $x' = y$. 
            Hence, $m_{v_0'} + m_{v_1'} \leq k-3$, and it reduces to~\ref{itm:leq_k-3}. 
            So we may assume $a_y \leq a-1$. 
            If $a_{v_0} = a+1$, then $a_{v_0} + a_y + a_{v-1} \leq 2(a+1) + a-1 \leq 3a+1$. 
            This implies that there exists $v' \in \Omega \setminus \set{v_0, y, v_1}$ such that $a_{v'} \geq \ceil{\frac{ak+2 - 3a-1}{k-3}} = a+1$. 
            Take another interior vertex $v'' \in \Omega \setminus \set{v_0, y, v_1, v'}$. 
            Let $v_0' = v_0, v_1' = v'$, and $x' = v''$. 
            Then $m_{v_0'} + m_{v_1'} \leq k-3$, and it reduces to~\ref{itm:leq_k-3}. 
            So we may assume $a_{v_0} \geq a+2$. 
            Take two distinct interior vertices $v', v'' \in \Omega \setminus \set{v_0, y, v_1}$. 
            Let $v_0' = v_0, v_1' = v'$, and $x' = v''$. 
            Then $m_{v_0'} + m_{v_1'} \leq k-3$, and it reduces to~\ref{itm:leq_k-3}.

            If $x$ is not adjacent to $v_0, v_1$ or $y$, then we have a path $v_0 \sim y \sim v_1 \sim z \sim x$ or $x \sim z \sim v_0 \sim y \sim v_1$. 
            If $k \geq 7$, then there exists $v' \in \Omega \setminus \set{v_0, y, v_1, z, x}$ adjacent to $v_0$ or $v_1$. 
            By~\cref{lambda at most minimum degree}, we have $a_{v'} \geq a$. 
            Let $v_0' = v_0, v_1' = v_1$, and $x' = v'$. 
            Then $m_{v_0'} + m_{v_1'} \leq k-3$, and it reduces to~\ref{itm:leq_k-3}. 
            So we may assume $k = 5$. 
            By~\cref{lambda at most minimum degree}, we have $a_x \geq a$ and $a_y, a_z \geq a-1$. 
            If $a_y \geq a$, then we take $v_0' = v_0, v_1' = v_1$, and $x' = y$. 
            Then $m_{v_0'} + m_{v_1'} \leq 5-3$, and it reduces to~\ref{itm:leq_k-3}. 
            So we may assume $a_y = a - 1$. 
            Similarly, we may assume $a_z = a - 1$. 
            
            Suppose the path is $v_0 \sim y \sim v_1 \sim z \sim x$.
            If $a_{v_1} \leq a$, we define 
            $h_1(v) = 
            \begin{cases}
                1, & v \in \set{y, v_1, z}, \\
                0, & \textotherwise.
            \end{cases}
            $ 
            Then $\lambda_1(T) \leq R_T(h_1) \leq a < \sigma_4(a, k)$. 
            So we may assume that $a_{v_1} \geq a+1$. 
            If $a_{x} \leq a$, we define 
            $h_2(v) = 
            \begin{cases}
                1, & v \in \set{z, x}, \\
                0, & \textotherwise.
            \end{cases}
            $ 
            Then $\lambda_1(T) \leq R_T(h_2) \leq a < \sigma_4(a, k)$. 
            So we may assume that $a_{x} \geq a+1$. 
            We take $v_0' = v_0, v_1' = x$, and $x' = v_1$. 
            Then $m_{v_0'} + m_{v_1'} \leq 5-3$, and it reduces to~\ref{itm:leq_k-3}. 
            
            Suppose the path is $x \sim z \sim v_0 \sim y \sim v_1$. 
            If $a_{v_1} \leq a$, we define 
            $h_3(v) = 
            \begin{cases}
                1, & v \in \set{y, v_1}, \\
                0, & \textotherwise.
            \end{cases}
            $ 
            Then $\lambda_1(T) \leq R_T(h_3) \leq a < \sigma_4(a, k)$. 
            So we may assume that $a_{v_1} \geq a+1$. 
            Similarly, we have $a_{x} \geq a+1$. 
            We take $v_0' = x, v_1' = v_1$, and $x' = v_0$. 
            Then $m_{v_0'} + m_{v_1'} \leq 5-3$, and it reduces to~\ref{itm:leq_k-3}. 
    \end{enumerate}
    The proof is complete. 
\end{proof}

\section{Concluding remark}\label{section:Concluding remark}

In~\cref{thm:lower_bound_2}, we prove the Li--Yau type Dirichlet eigenvalue estimates on trees based on inscribed radius.
However, we cannot obtain such good estimates for the first non-zero Laplacian eigenvalue of trees.
Mohar~\cite{mohar1991EigenvaluesDiameterMean} proved that $\mu_2(G) \geq \frac{4}{nD}$ for a graph of order $n$ and diameter $D$.
He provided an example as follows, which implies that this bound is sharp up to a constant factor and this bound cannot be improved to the form of $\frac{c}{D^2}$ even for trees.
Let $P_{k,t}$ be the tree of diameter $D=t+2$ which is obtained from the path $P_{t+1}$ by attaching $k$ new vertices to each of its two end vertices.
Mohar showed that
\begin{align*}
    \mu_2(P_{k, t}) < \frac{4}{(D + 2k - 1)(D - 2 - \epsilon)}
\end{align*}
for any $\epsilon>0$ as soon as $\frac{k}{t}$ is large enough.

Friedman~\cite{friedman1993some} first considered Faber--Krahn type inequality on graphs.

\begin{problem}
    Characterize all graphs with boundary which have minimal first Dirichlet eigenvalue among all graphs in a particular graph class.
\end{problem}

Bıyıkoğlu and Leydold~\cite{biyikouglu2007faber} answers the Faber--Krahn problem in the class of trees with given number of interior vertices and leaves or in the more refined class of trees with given degree sequence. 
One may consider the Faber--Krahn problem in the class of trees with given diameter.

\begin{problem}
    Characterize all trees which have minimal first Dirichlet eigenvalue for given number of vertices and given diameter.
\end{problem}

One may also consider the opposite of Faber--Krahn property, namely a graph with boundary achieving largest first Dirichlet eigenvalue among all graphs in a particular graph class. 

\begin{problem}
    Characterize all graphs with boundary which have largest first Dirichlet eigenvalue among all graphs in a particular graph class.
\end{problem}

Note that~\cref{subsec:b=ak+1,subsec:b=ak+2} gives partial answers in the class of trees with given number of interior vertices and leaves.
The other cases are still open.

\section*{Acknowledgements}
Huiqiu LIN was supported by the National Natural Science Foundation of China (No. 12271162, No. 12326372), and Natural Science Foundation of Shanghai (No. 22ZR1416300 and No. 23JC1401500) and The Program for Professor of Special Appointment (Eastern Scholar) at Shanghai Institutions of Higher Learning (No. TP2022031).
Da ZHAO was supported in part by the National Natural Science Foundation of China (No. 12471324, No. 12501459, No. 12571353), and the Natural Science Foundation of Shanghai, Shanghai Sailing Program (No. 24YF2709000). 

\bibliographystyle{abbrv}
\bibliography{ref}

\end{document}